\documentclass[10pt]{article}
\usepackage{amssymb,amsmath,amsfonts,mathrsfs,upgreek,textgreek,bbold} 
\usepackage{mathabx}

\usepackage[nottoc]{tocbibind}

\usepackage{graphicx} 

\usepackage[colorlinks=true]{hyperref}
\usepackage{comment}

\renewcommand\comment[1]{{\iffalse #1 \fi}}

\newtheorem{theorem}{Theorem}[section]

\newtheorem{lemma}{Lemma}[section]
\newtheorem{proposition}{Proposition}[section]

\newtheorem{assumption}{Assumption}[section]

\newtheorem{remark}{Remark}[section]
\newtheorem{notation}{Notation}[section]
\newtheorem{definition}{Definition}[section]
\newcommand{\diam}{{\,\rm diam\,}}

\newcommand{\real}{ {\mathbb R}   }
\newcommand{\torus}{ {\mathbb T}    }
\newcommand{\integer}{ {\mathbb Z}   }

\newcommand{\complex}{ {\mathbb C}   }

\newcommand{\cM}{ {\mathcal M}   }

\renewcommand{\Im}{\, {\rm Im}\,}
\renewcommand{\Re}{\, {\rm Re}\,}

\newcommand{\eproof}{\qed}

\newcommand\beq[1]{ \begin{equation}\label{#1} }
\newcommand{\eeq}{ \end{equation} }
\newcommand{\beqno}{ \[ }
\newcommand{\eeqno}{ \] }
\newcommand\beqa[1]{ \begin{eqnarray} \label{#1}}
\newcommand{\eeqa}{ \end{eqnarray} }
\newcommand{\beqano}{ \begin{eqnarray*} }
\newcommand{\eeqano}{ \end{eqnarray*} }

\newcommand\dfn[1]{ \begin{definition}\label{#1} }
\newcommand\edfn{ \end{definition} }
\newcommand\ass[1]{ \begin{assumption}\label{#1} }
\newcommand\eass{ \end{assumption} }

\newcommand\notat[1]{ \begin{notation} \label{#1} 
 }
\newcommand\enotat{\end{notation}}

\newcommand\rem{\begin{remark} 
\rm 
}
\newcommand\erem{\end{remark} 
}

\newcommand{\proof}{\noindent{\bf Proof\ }}
\newcommand\equ[1]{{\rm (\ref{#1})}}
\newcommand{\nl}{{\smallskip\noindent}}

\newcommand{\Giu}{{\bigskip\noindent}}
\newcommand{\noi}{{\noindent}}
\newcommand{\qed}{\hskip.5truecm
\vrule width 1.7truemm height 3.5truemm depth 0.truemm
\par\Giu}
\newcommand{\qedeq}{\hskip.5truecm
\vrule width 1.7truemm height 3.5truemm depth 0.truemm}

\newcommand\bolla{{\tiny ${}^\bullet\,$}}

\newcommand\casitwo[4]{ \left\{  \begin{array}{ll}
 {#1} & \mbox{ {\rm if} ${#2}$} \\
 {#3} & \mbox{ {\rm if} ${#4}$}
 \end{array} \right.}

\newcommand{\e}{\varepsilon}

\newcommand{\sign}{{\rm sign}}
\renewcommand{\a }{\alpha }
\renewcommand{\b }{\beta }
\newcommand{\s }{\sigma }
\newcommand{\ii }{{\rm i} }
\renewcommand{\d }{\delta }

\newcommand{\g }{\gamma}
\newcommand{\f }{\varphi}

\renewcommand{\l }{\lambda }

\newcommand{\C}{\mathbb{C}}

\newcommand{\z }{z}

\newcommand\ttg{{\mathtt g}}

\newcommand\ttq{{\mathtt q}}
\newcommand\ttp{{\mathtt p}}

\def\R{\mathbb R}
\def\T{\mathbb T}

\def\dst{\displaystyle}
\def\bks{\, \backslash\, }
\def\meas{{\rm\, meas\, }}

\newcommand\modulo{|}
\newcommand\ham{\mathtt{H}}

\newcommand{\Hpend}{{\ham}}
\newcommand{\bHpend}{{\bar\ham}}

\newcommand\lalla{\upmu}
\newcommand\lella{\upmu_{\rm o}}
\newcommand{\xx}{{\theta}}
\newcommand\loge{\uplambda} 

\newcommand\blogemax{\bar\loge_{{}_{\rm max}}}

\newcommand\sa{\theta} 
\newcommand\Sa{\Theta} 

\newcommand\Gm{{\mathtt G}}
\newcommand\GO{{\bar{\mathtt G}}}

\newcommand{\ch}{{\hat \eta}}

\newcommand{\suca}{\upepsilon}
\newcommand\morse{\upbeta}
\newcommand{\Ro}{{\mathtt R}} 
\newcommand{\ro}{{\mathtt r}}
\newcommand{\so}{{\mathtt s}}

\newcommand\act{I}  
\newcommand\ang{\f}  

\newcommand\Fiq{\Upphi^i}   

\newcommand\cFiq{\check\Upphi^i}  
\newcommand\Bu{{\mathcal B}}

\newcommand{\acci}{a}

\newcommand{\logemax}{{\loge_{{}_{\rm max}}}}

\newcommand\cin{\upnu}

\newcommand{\ts}{\textstyle}

\usepackage{appendix}

\renewcommand\ln{\log}
\newcommand\st{{\rm \ s.t.\ }}

\newcommand\cc{{\bf c}}
\newcommand\hcc{{\hat \cc}}

\newcommand\ddi\updelta			

\newcommand\bfco{{\bf c_{{{}_0}}}}
\newcommand\bfcu{{\bf c}_{{{}_1}}}
\newcommand\bfcd{{\bf c}_{{{}_2}}}
\newcommand\bfct{{\bf c}_{{{}_3}}}
\newcommand\bfcq{{\bf c}_{{{}_4}}}
\newcommand\bfcc{{\bf c}_{{{}_5}}}
\newcommand\bfcs{{\bf c}_{{{}_6}}}

\newcommand\tgo{{\ttg_{\rm o}}}
\newcommand\xxx{{\tau}}

\makeatletter
\newcommand{\pushright}[1]{\ifmeasuring@#1\else\omit\hfill$\displaystyle#1$\fi\ignorespaces}
\newcommand{\pushleft}[1]{\ifmeasuring@#1\else\omit$\displaystyle#1$\hfill\fi\ignorespaces}
\makeatother

\title{\bf 
Complex Arnol'd--Liouville maps
}

\begin{document}

\date{May 11, 2023}

\author{ 
\footnotesize L. Biasco  \& L. Chierchia
\\ \footnotesize Dipartimento di Matematica e Fisica
\\ \footnotesize Universit\`a degli Studi Roma Tre
\\ \footnotesize Largo San Leonardo Murialdo 1 - 00146 Roma, Italy
\\ {\footnotesize luca.biasco@uniroma3.it, luigi.chierchia@uniroma3.it}
\\ 
}

\maketitle

\begin{abstract}\noindent\footnotesize
We discuss the holomorphic  properties of the complex  continuation of 
 the classical Arnol'd--Liouville action--angle variables for real analytic 1 degree--of--freedom Hamiltonian systems depending 
on external parameters in suitable `generic standard form', with particular regard to the behaviour near separatrices.  

\nl

\nl
MSC2010 numbers: 37J05, 37J35, 37J40, 70H05, 70H08, 70H15 
\\
Key words: Hamiltonian systems. Action--angle variables. Arnol'd--Liouville integrable systems. Complex extensions of symplectic variables.   KAM Theory.
\end{abstract}

\section{Introduction}

Vladimir Igorevi\v c Arnol'd, in his didactic masterpiece \cite{A89} where most of us learned modern Mechanics,  explained, in precise mathematical terms, how to construct action--angle variables for an integrable Hamiltonian  system with compact energy levels (for a different  construction, see, also, \cite{HZ}). In this paper we consider real--analytic, one degree--of--freedom Hamiltonian systems depending on adiabatic invariants and discuss the fine holomorphic properties of the complex Arnold--Liouville map, especially near their singularities.

\nl
Analytic properties of the action--angle map are difficult to be found in the mathematical literature  despite their   interest {\sl per se} and, especially, in view of their 
relevance in modern perturbation theory; see, however, A.~Neishtadt's remarkable Thesis \cite{Nei89} (in Russian).
Indeed, in real--analytic theories (such as averaging theory, KAM theory, Nekhoroshev like theorems, {\sl etc.}\footnote{Compare  \cite{AKN} and references therein for general information.}),  it is necessary to control complex domains, whose characteristics appear explicitly in iterative  perturbative constructions. This is even more 
relevant if one needs to have holomorphic   information {\sl arbitrarily} close to the singularities of the action--angle variables (separatrices). 

\nl
For example, in real--analytic models of Arnol'd diffusion, starting with Arnol'd' pioneering 1964 paper \cite{A64},
one often  considers perturbations of the simple pendulum\footnote{\label{AD}
For references on Arnol'd diffusion,  besides \cite{A64}; compare, also, \cite{CG}, \cite{Z}, \cite{Ma}, \cite{T}, \cite{BKZ}, \cite{DLS}, \cite{CdL22}, \cite{KZ} among many other interesting results.}, while the analysis presented here allows for a real--analytic generic class of one--degrees--of--freedom Hamiltonians, which is a significant generalization of pendulum like models. Furthermore, in connection with the formidable problem of Arnol'd
diffusion in generic real--analytic {\sl a--priori} stable systems, the fine complex analytic understanding of the integrable limit appears to be as an essential tool.

\nl
Another example concerns the detection of primary and secondary Lagrangian tori in generic analytic nearly--integrable systems in phase space regions very close or inside the separatrices arising near simple resonances. In fact,   besides the 
usual KAM primary tori (i.e., the tori which are deformation of integrable ones at a distance $\gg\! \sqrt\e$ from separatrices), 
one expects the appearance of more tori very close to separatrices and secondary tori of different homotopy inside separatrices; compare, \cite[\S~6.3.3--C]{AKN} and \cite{MNT}. In this kind of analysis it also essential to have {\sl uniform estimates} over the relevant analyticity parameters, which is one of the main issues addressed in this paper. 
For more information on this subject,
see \cite{LW} for lower dimensional tori and \cite{BCnonlin}, \cite{BClin2}, \cite{BCuni}, and \cite{BCsecondary}, where the authors develop a `singular KAM Theory' for real--analytic generic natural systems, which shows, in particular,  that the total Liouville measure of Lagrangian KAM tori has a density larger than $1-\e |\log\e|^c$. For such results, the analytic tools developed in this paper play an essential r\^ole.

\nl
The main results of this paper are Theorem~\ref{glicemiak} and Theorem~\ref{barbabarba}. \\
In Theorem~\ref{glicemiak} it is shown
that near separatrices the actions, regarded as functions of  the  energy, have a special `universal' representation 
in terms of analytic functions and of logarithms for  energies  arbitrarily close to their singular values   (namely, the energy of the separatrices). Such a representation allows to give several estimates on the (derivatives of the)  action  functions.
\\
In Theorem~\ref{barbabarba} we compute the analyticity radii of the action--angle variables in arbitrary neighborhoods of separatrices and see how they behave in terms of a (suitably rescaled) distance from separatrices.\\
In \S~\ref{convexity}
 we investigate the convexity of the energy functions (defined as the inverse  of the action functions) near separatrices, showing that, in particular cases (in the outer regions outside the main separatrix, and in the case the potential is `close' to a cosine),
 the convexity is strictly defined, while in general it can be shown that inside separatrices there are inflection points. 
 
\nl
For completeness we include in Appendix~\ref{oceania} a rather standard hyperbolic normal form near hyperbolic equilibria.
 Appendix~\ref{australia} contains the proofs of two simple lemmata.
 
 \nl
 Finally, le us point out that  some results of this paper are not new. Indeed, similar analytic estimates
were obtained by different methods by A.I. Neishtadt in \cite[Ch. 3, Sect 7]{Nei89}, which, however,  
does not contain explicit estimates of constants related to analyticity
properties as discussed here.

\section{1D Hamiltonians in standard form}
\label{preliminari}

Consider a 1 degree--of--freedom real analytic   Hamiltonian system possibly parameterized by external  parameters with Hamiltonian function   
\beq{sting}
{\rm H}({\rm p},{\rm q_1})=\mathrm h({\rm p})+\e
\mathrm f({\rm p},{\rm q}_1)\,, 
\qquad
({\rm p},{\rm q_1})\in
\mathrm D\times \torus^1\,,
\eeq
where $n\ge 1$, $\mathrm D\subset \R^n$ is a bounded domain, $\torus^n:=\real^n/(2\pi \integer^n)$, 
$\e\ge 0$ a `perturbation parameter', $\hat {\rm p}=({\rm p_2},...,{\rm p}_n)$ are  external `adiabatic invariants',  
$\mathrm h$ and $\mathrm f$ are real analytic functions.

\nl
An important example is when $n\ge 2$ and ${\rm H}$ is the `secular' Hamiltonian of a nearly--integrable system
$(\ham,{\rm D}\times \torus^n)$ with 
\beq{ham} 
\ham=h(y)+\e f(y,x)
\eeq
real--analytic, after averaging (for  small $\e$) around a simple resonance 
\beq{resonance}
{\mathcal R}_k:=\{y\in {\rm D}: \partial_y h(y)\cdot k=0\}
\eeq
for some non vanishing $k\in\integer^n$, and after a suitable linear symplectic change of variables such that ${\rm q}_1=k\cdot x$ is the resonant angle and such that 
$\partial_{\mathrm p_1}\mathrm h|_{\mathrm p_1=0}= 0$. Here, 
`simple resonance' means that in the fixed neighborhood of ${\mathcal R}_k$ where averaging is performed,  there are no other independent resonant relations
$  \partial_y h(y)\cdot \ell=0$ for some $\ell$ independent of $k$ (and of not too high order); `secular' means that ${\rm H}$ is obtained disregarding the high order perturbation obtained after averaging.

\rem
As just seen in the above example, the case when  $\partial_{\mathrm p_1}\mathrm h(\mathrm p^0)= 0$ for some point $\mathrm p^0\in {\rm D}$ appears naturally in perturbation theory. \\
On the other hand, the case when
 $\partial_{\mathrm p_1}\mathrm h(\mathrm p^0)\neq 0$ is trivial:
indeed, by the Implicit Function Theorem,
 for values of the energy $E$ close to 
$\mathrm h(\mathrm p^0)$, $\mathrm p$ close to
$\mathrm p^0$  and $\e$ small enough, 
 there exists a function $v(E,\hat{\mathrm p},\mathrm q_1)$
 such that $\mathrm  h(v(E,\hat{\mathrm p},\mathrm q_1),\hat{\mathrm p},
 \mathrm q_1)=E$; therefore one can  define new action variables
$$I_1:=\int_0^{2\pi}v(E,\hat{\mathrm p},\mathrm q_1)
 d \mathrm q_1\,,\qquad \hat I:=\hat{\mathrm p}\,,$$
which, by the classical Arnol'd--Liouville construction,  can be completed into a symplectic transformation $\phi: (I,\f)\to {\rm (p,q)}$ such that 
 $\mathrm H\circ \phi = h(I)$ with $h$ real--analytic; see \cite[Ch. 10]{A89} for details.
\erem
Next, we show that, in general,  the Hamiltonian ${\rm H}$ in \equ{sting}, in a neighborhood of a critical point $\mathrm p^0$ of $h$,   can be symplectically put 
into a `standard form', which generalizes the features of the standard pendulum and 
it is particular suited to study  its (complex) Arnol'd--Liouville action--angle variables. The precise quantitative definition 
of `Hamiltonian in standard form' is given in the following two definitions.

\begin{definition}\label{buda}
A $C^2(\T,\R)$ Morse function $F$ with distinct critical values  is called  {\bf $\b$--Morse},  with $\b>0$,  if 
\beq{ladispoli}
\min_{\sa\in\T} \ \big( |F'(\sa)|+|F''(\sa)|\big) 
\geq\b \,,
\quad
\min_{i\neq j } 
|F(\sa_i)-F(\sa_j)|
\geq \b \,, 
\eeq
where $\sa_i\in\T$ are the  critical points of $F$.
\end{definition}
To formulate the next definition we need some notation:
Given $D\subseteq \real^m$, and $r>0$, we denote by $D_r$  the complex neighborhood of $D$ given by
$$D_r := \bigcup_{z\in D} \{y\in\complex^m\st |y-z|<r\}\,,$$
and, for $s>0$,  by $\torus^m_s$  the complex neighborhood of width $2s$ of $\torus^m$ given by
\beq{muu}
\torus^m_s:=\{x=(x_1,...,x_m)\in\complex^m:\ \ |\Im x_j|<s\}/(2\pi \integer^m)
\,.
\eeq

\begin{definition}\label{morso}
Let $\hat D \subseteq \R^{n-1}$ be a bounded domain,  $\Ro>0$ and $D:=  (-\Ro ,\Ro ) \times\hat  D
$. We say that the real analytic Hamiltonian $\Hpend$ is in Generic Standard Form (in short, `standard form')
with  respect to standard symplectic variables $(p_1,q_1)\in (-\Ro ,\Ro )\times\torus$ and  `external actions' 
$\hat p=(p_2,...,p_n)\in \hat D$ if $\Hpend$ has the form
 \beq{pasqua}
  \Hpend(p,q_1)=
\big(1+ \cin(p,q_1)\big) p_1^2  
  +\Gm(\hat p, q_1)
  \,,
\eeq
where $p=(p_1,\hat p)=(p_1,...,p_n)$ and:

\begin{itemize}

\item[\bolla]  $\cin$ and $ \Gm$ are real analytic functions defined on, respectively, $D_\ro\times\T_\so$ and $\hat D_\ro\times \T_\so$ for some $0<\ro\leq\Ro$ and $\so>0$;

\item[\bolla]  
$\Gm$ has zero average and there exists a 
 function $\GO$ (the `reference potential')  depending only on $q_1$ such that, for some $\morse>0$, 
 \begin{equation}\label{A2bis}
\GO\  \  \mbox{is} \ \
\morse {\rm \text{--}Morse}\,,\qquad \langle \GO\rangle=0\,;
\end{equation}
\item[\bolla] 
 the following estimates hold:
 \beq{cimabue}
 \left\{\begin{array}{l} \dst \sup_{\torus^1_\so}|\GO|\le \suca\,,\\
 \dst \sup_{\hat D_\ro\times \torus^1_\so}|\Gm-\GO| \leq
\suca
 \lalla
\,,\quad{\rm for\ some}\quad 0<\suca\le \ro^2/2^{16} 
\,,\ \ 0\le \lalla<1\,,
\\
\dst \sup_{D_\ro\times \torus^1_\so}|\cin| \leq
\lalla\,.
\end{array}\right.
\eeq
\end{itemize}
\end{definition}
We shall call $(\hat  D,\Ro,\ro,\so,\morse,\suca,\lalla)$ {\sl the analyticity characteristics of $\Hpend$
with respect to the unperturbed potential $\GO$}.
\rem\label{trivia}  (i) The Hamiltonian in standard form $\ham$ retains the basic features of the standard pendulum or more precisely of a natural system 
with a generic periodic potential, having,  in particular all equilibria  on the $p_1=0$ axis in the $(p_1,q_1)$--phase space.
%
%
%

\nl
(ii)
 If $\Hpend$ is in standard form, then the parameters $\morse$ and $\suca$ satisfy the relation\footnote{By 
 \equ{ladispoli} and \equ{cimabue}, $\morse \le | \GO(\sa_i)-  \GO(\sa_j)|\le 2 \max_\torus| \GO|\le 2\suca$.
 }
\beq{sucamorse}
\suca/\morse\ge 1/2\,.
\eeq 
Furthermore,  one can always
fix  $\upkappa\geq 4$  such that:
  \begin{equation}\label{alce}
1/\upkappa\leq \so\leq 1\,,\qquad
1\leq
\Ro/\ro\leq \upkappa\,,\qquad
1/2\leq
 \suca/\morse 
\leq \upkappa \,.
\end{equation}
Such parameter $\upkappa$ rules the scaling properties of these Hamiltonians and is the only constant (besides the dimension $n$)
on which the various constants depend. 

\nl
(iii)
The  critical points of a Morse function on $\torus$ (i.e., a function which has only non--degenerate critical points), by compactness, cannot accumulate, hence, they  are in a finite, even   number (alternately, a  relative maximum and a relative minimum). For $\b$--Morse functions one can easily estimate the number of critical points:

\nl
{\sl If $G$ is a $\b$--Morse function, then the number $2N$ of its critical points  does not exceed 
$\pi\sqrt{2\max_\R |G''|/\b}$}

\nl
\proof
If $\sa_i$ and $\sa_j$ are different critical points of $G$, then, by Taylor expansion at order two and by \equ{ladispoli}  one has
$\b\le|G(\sa_i)-G(\sa_j)|\le \frac12 (\max_\R|G''|) |\sa_i-\sa_j|^2$
which implies that 
\beq{dicembre}
\min_{i\neq j}|\sa_i-\sa_j|\ge\sqrt{2\b/\max_\R |G''|} 
\,,
\eeq
from which the claim follows at once. \qed
(iv) Of course, the constant $1/2^{16}$ appearing in the definition is quite arbitrary (as long as it is $\ll 1$).

\nl
(v) Hamiltonians in standard form have been  
investigated in \cite{BCnonlin}, \cite{BCuni} and \cite{BCsecondary}.
\erem

\begin{proposition}\label{rossana}
Let   ${\rm H}$ in \equ{sting} be a real analytic function and  assume that  at $\rm p^0\in D$  $\rm p_1\to h$ has a non--degenerate critical point\footnote{Explicitly: $\partial_{\mathrm p_1}\mathrm h(\mathrm p^0)=0$ and $\partial^2_{\mathrm p_1^2}\mathrm h(\mathrm p^0)\neq 0$.}\!. Assume also  that  
$\mathrm q_1\to\mathrm f(\mathrm p^0,\mathrm q_1)$ is a Morse function with distinct critical values.
Then,  for $\e$ small enough, ${\rm H}$ is symplectically conjugated to a Hamiltonian in standard form in a ($\e$--independent) neighborhood of $\{\mathrm p^0\}\times \torus^n$.
\end{proposition}

\proof 
Assume that $\mathrm h$ and $\mathrm f$ have holomorphic extension on, respectively, 
$\mathrm D_{\mathrm r}$ and
$\mathrm D_{\mathrm r}\times \T^1_{\mathrm s}$
for some $\mathrm r, \mathrm s>0$, and that 
$|\mathrm h|$ and $|\mathrm f|$ are uniformly 
bounded on their complex domains by some constant
$\mathrm M>0$.
Let us consider  $\mathrm H$ as a 1--degree--of--freedom  Hamiltonian 
in action--angle variables 
$(\mathrm p_1,\mathrm q_1)$,
depending on  parameters $\hat{\rm p}=(\hat{\rm p}_2,\ldots,\hat{\rm p}_n)$. 
\\
By assumption, there exist 
$\d,\b>0$ such that
$|\partial^2_{\mathrm p_1^2}\mathrm h(\mathrm p^0)| =\d$, and $F(\theta)=\mathrm f(\mathrm p^0,\theta)$ verifies \equ{ladispoli}.
By the Implicit Function Theorem\footnote{See, e.g., \cite[Appendix A]{Ch}.},
for
$c=c(n)>1$ large enough, setting
 $$\rho:=\d \mathrm r^3/c \mathrm M\leq \mathrm r/4\,,$$
 there exists a function
 $u(\hat{\mathrm p})$ holomorphic 
in
  $|\hat{\mathrm p}-\hat{\mathrm p}^0|<\rho$, 
  such that
 $\mathrm p_1^0=u(\hat{\mathrm p}^0)$,
 \beq{solero}
 \partial_{\mathrm p_1}\mathrm h(u(\hat{\mathrm p}),\hat{\mathrm p})= 0\qquad
 \mbox{and}\qquad
 |u(\hat{\mathrm p})-\mathrm p_1^0|\leq \rho\,,\ \ \forall\, 
 |\hat{\mathrm p}-\hat{\mathrm p}^0|<\rho\,.
 \eeq
Under the symplectic change of variables $(p_1,q_1)\to (\mathrm p_1,\mathrm q_1)$ given by 
$\mathrm p_1=u(\hat p)+p_1$,
$\mathrm q_1=q_1$ (with $\hat{\mathrm p}=\hat p$)
the  Hamiltonian $\rm H$ becomes
$$
H(p,q_1):=h(p)+\e f(p,q_1)\,,
$$
where
$h(p):=\mathrm h(u(\hat p)+p_1,\hat p)
$ and
$f(p,q_1):=\mathrm f(u(\hat p)+p_1,\hat p,q_1)$.
Note that  the function $h$ and $f$
are holomorphic and uniformly bounded
(in modulus) by $\mathrm M$
on, respectively, 
$\{|p_1|<\rho\}\times\{|\hat p-\hat {\mathrm p}^0|<\rho\}$
and
$\{|p_1|<\rho\}\times\{|\hat p-\hat {\mathrm p}^0|<\rho\}\times\T^1_{\mathrm s}$.
By \equ{solero}
we have that
\beq{solero2}
\partial_{p_1}
h(0,\hat p)=0\qquad
\forall\, 
|\hat p-\hat{\mathrm p}^0|<\rho\,.
\eeq
Moreover, by Cauchy estimates (and taking $c$ large enough)
we also have
$$
|\partial^2_{p^{2}_1}h(p)|\geq \d/2\,,\qquad
\forall\ 
|p_1|,\,|\hat p-\hat {\mathrm p}^0|<\rho
\,.
$$
Now we want to solve the equation
$\partial_{p_1}H=0.$
By \equ{solero2} we have 
$$\partial_{p_1}H(0,\hat p,q_1)|_{\e=0}=\partial_{p_1}
h(0,\hat p)=0\,.$$
By the Implicit Function Theorem
there exists
$c_*=c_*(n)\geq c>1$ large enough such that,
if $\sqrt\e\leq\d \mathrm r^2/c_* \mathrm M$,
then there exists a  function
$v(\hat p,q_1)$ holomorphic
on $$\{|\hat p-\hat {\mathrm p}^0|<\rho\}\times\T^1_{\mathrm s}$$
with $|v|\leq c_*\e \mathrm M/\d \mathrm r\leq \rho/4$, satisfying 
$\partial_{p_1}H(v(\hat p,q_1),\hat p,q_1)=0$.
Then let us perform the symplectic transformation
$p_1=v(\hat p,q_1)+p_1$, $q_1=q_1$
(with $\hat p=\hat {\rm p}$).
The new Hamiltonian $H(v(\hat p,q_1)+p_1,\hat p,q_1)$
 is holomorphic on 
 $$\{|p_1|<\rho/4\}\times
 \{|\hat p-\hat{\mathrm p}^0|<\rho\}\times\T^1_{\mathrm s}\,,$$
and is given by\footnote{For brevity we write $u$ and $v$
 instead of $u(\hat p)$ and  $v(\hat p,q_1)$,
 respectively.}
\beqa{gingerina}
H(v+p_1,\hat p,q_1)
&=&
H(v,\hat p,q_1)
+p_1^2 \int_0^1 (1-t)
\partial^2_{p^{2}_1}H(v+tp_1,\hat p,q_1)\, dt
\nonumber
\\
&=&
g(\hat p)+g_0  \ham\,,\qquad \ham:=
\big(1+ \cin(p,q_1)\big) p_1^2  
  +\Gm(\hat p, q_1)\,,
\eeqa
where
$g(\hat p):=h(0,\hat p)$,
 $g_0:=\frac12\partial^2_{\mathrm p_1^2}\mathrm h(\mathrm p^0)$
 and
\beqa{ginger}
\cin(p,q_1)\!\! \!&:=&\!\!\!\frac1{g_0}
\int_0^1 (1-t)\big(\partial^2_{\mathrm p_1^2}\mathrm h
(u+v+t p_1,\hat p)
-\partial^2_{\mathrm p_1^2}\mathrm h(\mathrm p^0)
+\e \partial^2_{\mathrm p_1^2}\mathrm f(u+v+tp_1,\hat p,q_1)
\big)dt,
\nonumber
\\
\Gm(\hat p, q_1) \!\! \!&:=&\!\!\!\frac1{g_0}
\int_0^1 (1-t)\partial^2_{\mathrm p_1^2}\mathrm h
(u+tv,\hat p)v^2\,dt
\, +\, 
\e \mathrm f(u+v,\hat p,q_1)\,,
\nonumber
\\
\GO(q_1) \!\! \!&:=&\!\!\!\frac\e{g_0}\ 
\mathrm f(\mathrm p^0, q_1)\,.
\eeqa
By \equ{solero},\equ{solero2}, Cauchy estimates,
the facts that 
$$\sqrt\e\leq\d \mathrm r^2/c_* \mathrm M\,,\quad
4|v|\leq\rho=\d \mathrm r^3/c \mathrm M\,,$$
 and $|u+v+tp_1-\mathrm p_1^0|\leq 2\rho\leq\mathrm r/2$ for every
 $0\leq t\leq 1$,
and noting that $|g_0|=\d/2$
and $\d\leq 2\mathrm M/\mathrm r^2$ (by Cauchy estimates),
it is easy to see that the Hamiltonian $\ham$ 
in \equ{gingerina}--\equ{ginger}
is in standard form according to Definition~\ref{morso}
with analyticity characteristics:
\beqno\ts
\hat D:=\{|\hat p-{\mathrm p}^0|<\frac{\rho}8\}\,,
\
\Ro=\ro:=\frac{\rho}8\,,
\
\so:=\min\{s,1\}\,,
\
\morse:=
\frac{2\e \b}{\d}\,,
\
\suca:=\frac{3\e \mathrm M}{\d}\,,
\
\lalla:=\frac{144}{c}
\,,
\eeqno
for a suitable constant  $c>144$.
In particular condition $0<\suca< \ro^2/2^{16}$ in \equ{cimabue}
is satisfied taking
$$
\e\leq\frac{\d^3 \mathrm r^6}{2^{24}c^2\mathrm M^3}\,.
$$
Finally, taking $\upkappa:=\max\{4,1/s\}$, \equ{alce}
holds. \qed


\section{Analytic properties  of actions at critical energies}

In the rest of the paper we shall investigate the complex analytic properties of the action--angle variables for a Hamiltonian in standard form. 

\nl
In this section 
we show that near separatrix the actions regarded as functions of  the  energy $E$   have a quite special `universal' representation 
(in terms of analytic functions and of logarithms)
for  energies  close to their singular values   (namely, the energy of the separatrices).

\nl\nl
Let   $\Hpend$ be a  Hamiltonian in standard form  (Definition~\ref{morso}), 
let $\bar \sa_0$ be  the unique absolute maximum of the reference potential $\GO$ in $[0,2\pi)$, and let $2N$ be the number of its critical points (compare  Remark~\ref{trivia}--(iii)). 
Then, the relative strict non--degenerate maximum and minimum
 points of $\GO$, $\bar\sa_i\in[\bar \sa_0,\bar \sa_0+2\pi]$, ($0\leq i\leq 2N$)
follow in alternating order,
$\bar\sa_0< \bar\sa_1<\bar\sa_2<\ldots <\bar\sa_{2N}:=\bar\sa_0+2\pi$,
 in particular,  {\sl $\bar\sa_i$ are
 relative maxima/minima points 
 for $i$ even/odd}. \\
Since $\GO$ is a $\morse$--Morse  function, 
 the corresponding  critical energies are distinct; let us denote them  by $\bar E_i:=\GO(\bar\sa_i)$. Hence,  
{\sl $\bar E_{2N}=\bar E_0$ being the unique  global maximum}.

\nl
By the Implicit Function Theorem,
 for $\lalla$ small enough,
 we can continue  the $2N$ critical points $\bar\sa_i$ of $\GO$    obtaining $2N$ 
critical points
$\sa_i=\sa_i(\hat{p})$
of $\Gm(\hat{p},\cdot)$ for $\hat{p}\in\hat D$;
the corresponding distinct  critical energies become
 \begin{equation}\label{faciolata2}
E_j(\hat p):=\Gm(\hat p, \sa_j(\hat p))\,.
\end{equation}
In fact, the following simple  lemma -- proven in Appendix~\ref{cento} -- based on the Implicit Function Theorem holds.

\begin{lemma}\label{enricone} 
If \footnote{$\upkappa$ as in \equ{alce}.}
$\lalla\leq 1/(2\upkappa)^6$
then
the functions   $\sa_i(\hat{p})$
and $E_i(\hat{p})$ are  real analytic in $\hat{p}\in \hat D_\ro$
and
\begin{equation}\label{october}\ts
\sup_{\hat{p}\in \hat D_\ro}|\sa_i(\hat{p})-\bar\sa_i| \leq
\frac{2\suca\lalla}{\morse \so}\,,\qquad 
\sup_{\hat{p}\in \hat D_\ro} |E_i(\hat{p})-\bar E_i|
 \leq 3\upkappa^3 \suca\lalla\,.
\end{equation}
Furthermore, the relative order of  $\sa_i(\hat{p})$ and $E_i(\hat{p})$ is, for every $\hat{p}\in \hat D_\ro$, the same as that of, 
respectively, $\bar\sa_i$ and $\bar E_i$.
\end{lemma}

\nl
Now, let  $\hat p\in\hat D$ and consider the following  phase space 
of the 1D Hamiltonian system governed by $\Hpend$: 
\beq{insalatina}
\cM=\cM(\hat p):=
\{
(p_1,q_1)\in (-\Ro -\ro,\Ro+\ro)\times\T\ \ \mbox{s.t.}\ \ 
\Hpend(p_1,\hat p,q_1)<\Ro^2+\Ro\ro
\}\,.
\eeq
Then, $\cM$ 
decomposes in  $2N+1$ 
open connected components $\cM^i=\cM^i(\hat {p})$, with  $0\le i\le 2N$, plus a zero measure singular   set $S=S(\hat p)$
 formed by  the $2N$ connected
separatrices and the $2N$ critical points:
\beq{enkidu} 
\cM=\cM(\hat p)=\bigcup_{i=0}^{2N} \cM^i  \ \cup\ S =
\bigcup_{i=0}^{2N} \cM^i(\hat {p}) \ \cup\ S(\hat p)\,.
\eeq
 Although the sets $\cM$ and $S$ depend upon the dumb actions $\hat p$, for $\lalla$ as in 
Lemma~\ref{enricone},
one sees easily that  
\beqa{gonatak}
\big(-\Ro -\frac\ro3,\Ro+\frac\ro3\big)\times\T
\subseteq\cM\subseteq     
\big(-\Ro -\frac\ro2,\Ro+\frac\ro2 \big)\times\T\,.
\eeqa
For the labelling of the domains
 $\cM^i(\hat p)$ we shall adopt the following conventions:
 
 \nl
$\cM^{0}(\hat p)$ is the region below the lowest separatrix
and 
$\cM^{2N}(\hat p)$ is the region above the highest separatrix.
\\
For $1\leq i\leq 2N-1$ odd the closure of 
$\cM^{i}(\hat p)$ contains the minimum point 
$(0,\sa_i)$, while, for $2\leq i\leq 2N-2$ even
the boundary of
$\cM^{i}(\hat p)$
is formed by two connected components the inner one
 contains the maximum point 
$(0,\sa_i)$.

\nl
For $0\le i\le 2N$, define
\begin{eqnarray}
\label{mortazza}
&&E^{(i)}_-(\hat p):=E_i(\hat p)\,,\ \ 0\leq i\leq 2N\,,\quad
E^{(0)}_+(\hat p)=E^{(2N)}_+(\hat p):=\Ro^2+\Ro\ro\,,
\nonumber
\\
&&E^{(2j-1)}_+(\hat p):=\min\{ E_{2j-2}(\hat p),
E_{2j}(\hat p)\}\,,\ \ 
1\leq j\leq N\,,
\nonumber
\\
&&E^{(2j)}_+(\hat p):=\min\{ E_{2j_-}(\hat p), 
E_{2j_+}(\hat p)\}\,,\ \ 
1\leq j< N\,,
\end{eqnarray}
 with 
$\dst 
j_-:=\max\{ i<j\ \ {\rm s.t.}\ \ \bar E_{2i}>
 \bar E_{2j} \}$, 
$j_+:=\min\{ i>j\ \ {\rm s.t.}\ \ \bar E_{2i}>
 \bar E_{2j} \}$. 

\nl
Then, for $\hat p\in\hat D$ fixed,  and for every $0\le i \le 2N$, we can define {\sl the action functions}
\beq{azione}
E\in (E^{(i)}_-(\hat p),E^{(i)}_+(\hat p)) \to\act_1^{(i)}(E,\hat p)
\eeq
by the standard 
Arnol'd--Liouville's  formula
\beq{4-1}
\act_1^{(i)}(E,\hat p):= \frac1{2\pi}\oint_{\g_i} \ p_1dq_1\,, 
\eeq 
where
\beq{poirot}
\g_i=\g_i(E;\hat p):=
\Hpend^{-1}(E;\hat p)
\cap \cM^i(\hat p)
\eeq
 is the smooth closed curve
in the plane $(q_1,p_1)$ with clockwise orientation\footnote{For the external curves ($i=0,2N$) the orientation is to the right
in $\cM^{2N}(\hat p)$, to the left in 
$\cM^0(\hat p)$.}.

\nl
Finally, we denote by $\bar \act_1^{(i)}(E)$
the action variables of the `unperturbed' Hamiltonian 
\beq{vespasiano}
\bHpend:=\Hpend|_{\lalla=0}:=p_1^2+\GO(q_1)\,,
\eeq
and observe  that $\act_1^{(i)}(E,\hat p)$ 
 reduces to $\bar \act_1^{(i)}(E)$
 for $\lalla=0.$

\nl
The  main result of this section is the following

\begin{theorem}\label{glicemiak} Let $\Hpend$ be a   Hamiltonian in standard form as in Definition~\ref{morso},
let $\upkappa\ge 4$ be such that \equ{alce} holds and let $2N$ be the number of critical points of the reference potential $\GO$.
Then, there exists   a suitable  constant
$\cc=\cc(n,\upkappa)\geq  2^8\upkappa^3$ 
such that, if
 \begin{equation}\label{caviale2}
\lalla\leq 1/\cc^2\,,
\end{equation}   
then, for all $0\le i\le 2N$ and $\hat \act\in\hat D$, the action functions in \equ{azione}  verify the following properties.

\nl
{\rm (i)} {\bf(Universal behaviour at critical energies)} 
There exist  functions $\phi^{i}_-(\z,\hat \act),$ 
$\psi^{i}_-(\z,\hat \act)$ for  $0\le i\le 2N$, 
and,  functions $\phi^{i}_+(\z,\hat \act),$ 
$\psi^{i}_+(\z,\hat \act)$, for $0<i<2N$,
 which are 
real analytic in 
$\{\z\in\complex: |\z|< 1/\cc\}\times  \hat D_{\ro/2}$  and satisfy
\begin{equation}\label{LEGOk}
\act_1^{i}\big(E_\mp^{i}(\hat \act)\pm \suca \z, \,\hat \act\big)=
\phi^{i}_\mp(\z,\hat \act) +\psi^{i}_\mp(\z,\hat \act)\ \z \log \z\, 
 \,,\quad \forall \ 0<\z<{ 1/\cc}\,,\,
 \hat \act\in \hat D\,.
\end{equation}
On  $\{\z\in\complex: |\z|< 1/\cc\}\times  \hat D_{\ro/2}$
the functions $\phi^{i}_\pm(\z,\hat \act)$,
$\psi^{i}_\pm(\z,\hat \act)$ satisfy:
\beq{pappagallok}
\begin{array}{l}
\dst
\sup_{|\z|<1/\cc, \,\hat \act\in \hat D_{\ro/2}}\big(
|\phi^{i}_\pm|+|\psi^{i}_\pm|\big)
\leq \cc
\sqrt\suca 
\,,
\\
\dst
\sup_{|\z|<1/\cc, \,\hat \act\in \hat D_{\ro/4}}
\big(|\partial_{\hat \act}\phi^{i}_\pm|
+|\partial_{\hat \act}\psi^{i}_\pm|\big)
\leq \cc
\lella
\,,\qquad \lella:=
{\ts \frac{\sqrt\suca}{\ro}\lalla} \stackrel{\eqref{cimabue}}\leq 
2^{-8}\lalla
\,.
\end{array}
\eeq
Moreover,
 \begin{equation}\label{trippa}
 | \phi^{i}_\pm -  \bar\phi^{i}_\pm|\,,\ 
| \psi^{i}_\pm -  \bar\psi^{i}_\pm|\ \leq\ 
\cc\sqrt\suca \lalla\,,
\end{equation}
where $\bar\phi^{i}_\pm:=\phi^{i}_\pm|_{_{\lalla=0}}$
and 
$\bar\psi^{i}_\pm:=\psi^{i}_\pm|_{_{\lalla=0}}$.

\nl
{\rm (ii)} {\bf (Limiting critical values)}
The following bounds at the limiting critical energy values hold:
\beq{ciofecak}
 \begin{array}{l}
|\psi^{i}_+(0,\hat \act)|\geq \sqrt\suca/ \cc \,,
 \quad \phantom{r} 0< i< 2N\,,\quad  \forall \ \hat \act\in \hat D_{\ro/2}\,,\\
 |\psi^{2j}_-(0,\hat \act)|\geq \sqrt\suca/\cc\,,
 \quad \ 0\leq j\leq N\,, \quad\   \forall  \ \hat \act\in \hat D_{\ro/2}\,,
 \\
 \psi^{i}_+(0,\hat \act)>0 \,,
 \quad \phantom{r} 0< i< 2N\,,\quad  \forall \ \hat \act\in \hat D\,,\\
 \psi^{2j}_-(0,\hat \act)<0\,,
 \quad \ 0\leq j\leq N\,, \quad\   \forall  \ \hat \act\in \hat D\,,
 \end{array}
\eeq
while, in the case of relative minimal critical energies, one has, $\forall$ $\hat \act\in \hat D$, $0<\z<{ 1/\cc}$,
\begin{equation}\label{lamponek}
\phi^{2j-1}_-  (0,\hat \act)=0\,,
\qquad     
\psi^{2j-1}_-(\z,\hat \act)=0\,,
\qquad  \forall\ 1\le j\le N\,.
\end{equation}
\nl
{\rm (iii)}
{\bf (Estimates on  derivatives of  actions on real domains)}
The derivatives of the action functions on real domains satisfy the following estimates:
\begin{equation}\label{vana} 
\inf_{( E^{i}_-, E^{i}_+)}
\partial_E  \act_1^{i} \ge \frac{1}{\cc\sqrt\suca}\,,\qquad \forall\ \hat \act \in \hat D\,,\ \forall\ 0<i<2 N
\,;
\end{equation}

\begin{equation}\label{moldavater}
\min\big\{\partial_E  \act_1^{2N}\,,
\
\partial_E  \act_1^{0}\big\}
\ge \frac{1}{\cc\sqrt{E+\suca} }
\,,
\ \ 
\forall\,
E> E_{2N}\,,\  \forall\ \hat \act\in\hat D\,.
\end{equation}
{\rm (iv)}
{\bf  (Estimates on  derivatives of  actions on complex domains and perturbative bounds)}
For $\loge>0$ satisfying
 \begin{equation}\label{bassoraTH}
  \cc\lalla
  \leq
  \loge
  \leq 
  1/\cc
  \,,
\end{equation}
define the following complex energy--domains:
\beq{autunno2}
{\mathcal E}^{i}_\loge :=
\left\{
\begin{array}{ll}
 \{E\in\complex:\bar E^i_- -{\ts \suca/\cc}\ \,<\Re  E<\bar E^i_+ - \loge\suca\,,\   |\Im E|<{\ts \suca/\cc} \}
\,,
&i\ {\rm odd}\,,
\\
 \{E\in\complex:\bar E^i_- + \loge\suca <\Re  E<\bar E^i_+ - \loge\suca\,,\  |\Im E|<{\ts \suca/\cc}\}
\,,
&0,2N\neq i\ {\rm even}
\,,
\\
 \{E\in\complex:\bar E^i_- + \loge\suca <\Re  E
 <\bar E^i_+
 \,,\  |\Im E|<{\ts \suca/\cc} \}
\,,
& i=0,2N
\,.
\end{array}\right.
\eeq
Then, for $0\leq i\leq 2N$, the functions 
$
 \act_1^{i}$
and
$
 \bar \act_1^{i}$ are holomorphic 
on the domains
 $\mathcal E^{i}_\loge
\times\hat D_{\ro}$,  and satisfy the following estimates:
 \begin{equation}
\sup_{\mathcal E^{i}_\loge
\times\hat D_{\ro/4}}
|\partial_{\hat \act} \act_1^{i}|
\leq \cc^2\, 
\lella
\,,
\ \ 
\sup_{\mathcal E^{i}_\loge}
\big|
\partial_E \bar \act_1^{i}
\big|
\leq \cc^2\, 
\frac{|\log\loge|}{\sqrt\suca }
\,,\ \ 
\sup_{\mathcal E^{i}_\loge
\times\hat D_{\ro/2}}
\big|
\partial_E \act_1^{i}
-
\partial_E \bar \act_1^{i}
\big|
\leq 
\frac{\cc^2\lalla}{\loge \sqrt\suca} \,
\label{rosettaTH}\,.
\end{equation}
\end{theorem}

\rem\label{caciottella}

\noi
(i)  
Statements similar to 
\eqref{LEGOk} have bee also discussed in \cite[Lemma 7.2]{Nei89}, \cite{Nei}, and 
 \cite[Eq. (5.8)]{BFS}.  Analyticity at elliptic equilibria (see Eq. \equ{lamponek} above), was proven also in 
 \cite[Lemma 7.1]{Nei89}.

\nl
(ii) 
Condition \equ{caviale2}
implies the hypothesis of Lemma \ref{enricone}. 
\erem
In the rest of the paper we shall  use the following

\begin{notation}\label{lessdot}
Given $m,M\ge 0$,  we say that $m\lessdot M$ if there exists a constant $c=c(n,\upkappa)\geq 1$ such that 
$m\le cM$. We shall also say that a function $f$ is of order $M$, $f=O(M)$, if $|f|\lessdot M$ uniformly on its domain of definition. 
\end{notation}

\nl
\proof {\bf  of Theorem~\ref{glicemiak}}
For definiteness we consider  the case of $i=2j+1$ odd and, in particular, 
the case with
$E_{2j}(\hat p)<E_{2j+2}(\hat p)$.
The other cases can be treated in  the same way with the obvious changes.

\nl
Recalling \eqref{mortazza} we note that
\beq{mortazzona}
E_+(\hat p):=E_+^{(2j+1)}(\hat p)=E_{2j}(\hat p)\,,
\qquad
\bar E_+:=\bar E_+^{(2j+1)}=\bar E_{2j}\,.
\eeq
\nl
For every fixed $\hat p\in\hat D$ 
we denote by
\beq{palumbo}
E\in \big(E_{2j+1}(\hat p),E_{2j}(\hat p)\big)\to\Sa_{{}_{\!\star}}(E,\hat p)\,,
\ \mbox{resp.,}\quad
E\in \big(E_{2j+1}(\hat p),E_{2j+2}(\hat p)\big)\to\Sa^{{}^{\!\star}}(E,\hat p)\,,
\eeq
the (real analytic) inverse
of $\Gm(\hat p, q_1)$
on the interval $\big(\sa_{2j}(\hat p),\sa_{2j+1}(\hat p)\big)$, respectively
$\big(\sa_{2j+1}(\hat p),\sa_{2j+2}(\hat p)\big)$.
As usual,  a bar above functions means the limit  $\lalla=0$, namely
$E\in (\bar E_{2j+1},\bar E_{2j})\to\bar\Sa_{{}_{\!\star}}(E)$, respectively,
$E\in (\bar E_{2j+1},\bar E_{2j+2})\to\bar\Sa^{{}^{\!\star}}(E)$,
will denote
the (real analytic) inverse
of $\GO(q_1)$
on the interval $(\bar\sa_{2j},\bar\sa_{2j+1})$, respectively
$(\bar\sa_{2j+1},\bar\sa_{2j+2})$.
Then, the action function $\bar\act(E)=\bar\act^{(2j+1)}_1 (E)$
of $\bHpend$ in \equ{vespasiano} can be   written as
\beq{frascati}
\bar \act(E):=
\frac{1}{\pi}
\int_{\bar \Sa_{{}_{\!\star}}(E )}^{\bar \Sa^{{}^{\!\star}}(E )}
\sqrt{E-\GO(\sa)}\, d\sa\,,\qquad
E\in(\bar E_{2j+1},\bar E_{2j})\,,
\eeq
so that
\beq{frascati2}
\partial_E\bar \act(E):=
\frac{1}{2\pi}
\int_{\bar \Sa_{{}_{\!\star}}(E )}^{\bar \Sa^{{}^{\!\star}}(E )}
\frac{d\sa}{\sqrt{E-\GO(\sa)}} \,.
\eeq
\nl
We split the proof in four steps.

\nl
{\bf Step 1:} {\sl Explicit expression for the action functions}

\nl
In this first step we  will obtain an analogous of \eqref{frascati} for
$I(E)=\act^{(2j+1)}_1 (E)$ in \equ{4-1}, see
formula \eqref{musicalbox} below;
estimates \eqref{vana} will then follow easily.

\nl
Let us consider the
  equation
\begin{equation}\label{meaux}
 p_1=\frac{z}{\sqrt{1+ \cin(p,q_1 )}}\,.
\end{equation}
Note that by \equ{cimabue} and \equ{caviale2} we have
\begin{equation}\label{aterno}
\Re(1+ \cin(p,q_1 ))\geq \frac12\,,\qquad
\forall\ 
 (p,q_1)\in D_\ro\times \torus^1_\so\,,
\end{equation}
and, therefore, 
$\sqrt{1+ \cin(p,q_1 )}$ is well defined 
on $D_\ro\times \torus^1_\so$.

\begin{lemma}\label{amatriciana}
There exists a unique real analytic function
  $\tilde{\mathcal P}:(-\Ro,\Ro)_{\ro/4}\times\mathbb T_{\so}\times \hat D_{\ro}\, \to\, \C$
satisfying the bound
 \begin{equation}\label{3holes}
 |\tilde{\mathcal P}|_\dagger:=
\sup_{(-\Ro,\Ro)_{\ro/4}\times\mathbb T_{\so}\times \hat D_{\ro}} 
{\modulo}\tilde{\mathcal P}{\modulo}
\leq  2\lalla\Ro\leq \frac{\ro}{8}\,,
\end{equation} 
 and  such that
\begin{equation}\label{spiderman}
p_1=\mathcal P(z,q_1,\hat p ):=
z +\tilde{\mathcal P}(z,q_1,\hat p )
\end{equation}
solves \eqref{meaux}:
\begin{equation}\label{meauxter}
\mathcal P(z,q_1,\hat p )=
\frac{z}{\sqrt{1+ \cin(\mathcal P(z,q_1,\hat p ), \hat p,q_1 )}}\,.
\end{equation}
Moreover,
\begin{equation}\label{piove}
\mathcal P: 
(-\Ro,\Ro)_{\ro/4}\times\mathbb T_{\so}\times \hat D_{\ro}
\to (-\Ro,\Ro)_{\ro/2}\,.
\end{equation}
\end{lemma}
\proof
We first note that 
if $\tilde{\mathcal P}$ satisfies the first inequality in \eqref{3holes},
then, by \eqref{alce} and \eqref{caviale2},   it follows that it also satisfies the second one.
Therefore,
if $z\in (-\Ro,\Ro)_{\ro/4}$,
then
  $z+\tilde{\mathcal P}\in (-\Ro,\Ro)_{\ro/2}$
  and \eqref{piove} holds.
Let $\mathtt B$ denote the closed ball  of functions 
$\tilde{\mathcal P}$ satisfying
\eqref{3holes} and let 
 $\tilde{\mathcal P}=
 \tilde{\mathcal P}(z,q_1,\hat p )$ be the solution of 
the fixed point equation
\begin{equation}\label{meauxbis}
\tilde{\mathcal P}=\Phi(\tilde{\mathcal P}):=
 \Big( \big(1+  \cin(z+\tilde{\mathcal P},\hat p,q_1 )\big)^{-\frac12}-1\Big)z
\end{equation}
By \eqref{cimabue}, \eqref{alce}
and \eqref{aterno}, it follows
$$
|\Phi(\tilde{\mathcal P})|_\dagger
\leq 
\lalla(\Ro+\ro/4)
 \leq 2\lalla\Ro\,,
$$
and, therefore,  $\Phi(\mathtt B)\subseteq
\mathtt B$. In fact, $\Phi$ is a contraction:
Omitting for brevity to write
$\hat p, q_1$ and setting
$\xx(t)=(1-t)\tilde{\mathcal P}' +t\tilde{\mathcal P}$,
 we get
$$
\cin(z+\tilde{\mathcal P})-
\cin(z+\tilde{\mathcal P}' )
=\big(\tilde{\mathcal P}-\tilde{\mathcal P}'\big)\int_0^1 \partial_{p_1} 
 \cin(z+\xx(t) )dt\,.
 $$
 Since $|\xx(t)|_\dagger\leq 2\lalla\Ro\leq \ro/{8}$
 and 
 $z+\xx(t)\in (-\Ro,\Ro)_{\ro/2}$
  for every $0\leq t\leq 1$,
   by \eqref{cimabue}
 and Cauchy estimates we get
 $|\partial_{p_1} 
 \cin(z+\xx(t) )|_\dagger\leq 2\lalla/\ro$
for any $0\leq t\leq 1$.
 Then, by \eqref{alce} and \eqref{caviale2},
$$
|\Phi(\tilde{\mathcal P})-\Phi(\tilde{\mathcal P}')|_\dagger
\leq 2 \big|\big(\cin(z+\tilde{\mathcal P},\hat p,q_1 )-
\cin(z+\tilde{\mathcal P}',\hat p,q_1 )\big)z\big|_\dagger
\leq \frac{8\lalla\Ro}{\ro}
|\tilde{\mathcal P}-\tilde{\mathcal P}'|_\dagger
\leq \frac12 
|\tilde{\mathcal P}-\tilde{\mathcal P}'|_\dagger
\,,
$$
and   \eqref{meauxbis} is solved  by the standard Contraction Lemma.
\eproof
Thus, for real values of $\hat p,$ $q_1,$ $E$
such that $0\leq E-\Gm(\hat p,q_1 )\leq \Ro+\ro/4$, we have that 
\begin{equation}\label{ummagamma}
p_1=\mathcal P\Big(\pm\sqrt{E-\Gm(\hat p,q_1 )},q_1,\hat p \Big)\quad
{\rm solves}  \quad
\Hpend (p_1,\hat p,                                                                                                                                                                                                                                                                                                                                                                                                                                                                                                                                                                                                                                                                                                                                                                                                                                                                                                                                                                       q_1 )=E\,,
\end{equation}
where the sign depends on whether $\pm  p_1\geq 0$.
By \eqref{3holes}, \eqref{spiderman} and Cauchy estimates,  for $z\in (-\Ro,\Ro)$, 
\begin{equation}\label{fragile}
\partial_z \mathcal P
\geq \frac12\,,
\end{equation}
so that for real $q_1,\hat p$, the real function  
$z\in (-\Ro,\Ro)\mapsto \mathcal P(z,q_1,\hat p)$ is  increasing.
Note also that
$\mathcal P(0,q_1,\hat p)= 0$.

\nl
Define the analytic function
\begin{equation}\label{4holes}
 \cin_\sharp(z,\sa,\hat p ):=
\frac{1}{2\sqrt{1+ \cin\big(\mathcal P(z,\sa,\hat p ),\hat p,\sa \big)}}
+
\frac{1}{2\sqrt{1+ \cin\big(\mathcal P(-z,\sa,\hat p ),\hat p,\sa \big)}}
-1\,.
\end{equation}
Notice that $\cin_\sharp$ is {\sl even} in $z$ and 
that\footnote{If $\Re (1+\cin)\geq 1/2$ (see \eqref{aterno}), then $|(1+\cin)^{-1/2}-1|\leq |\cin|.$},
by \eqref{piove} and   \eqref{cimabue},
\begin{equation}\label{rodimento}
\sup_{z\in (-\Ro,\Ro)_{\ro/4}}  {\modulo}\cin_\sharp(z,\hat p,\sa){\modulo}_{\hat D,\ro,\so}
\leq 
\sup_{D_\ro\times \torus^1_\so}|\cin| \leq
\lalla
\,.
\end{equation}
Then, by 
\eqref{rodimento}, \eqref{alce} and Cauchy estimates
we have
\begin{equation}\label{prurito}
\sup_{z\in (-\Ro,\Ro)_{\ro/8}}
|z\partial_z \cin_\sharp|_{\hat D,\ro,\so}
\leq
16\upkappa\lalla\,.
\end{equation}
We now need the following elementary 

\begin{lemma}\label{geronimo} Let $g:(-r,r)\to\real$ be an even function with  
holomorphic extension on 
 $[0,R]_r$. Then,
 one can define $G$ holomorphic on  $[0,R^2]_{r^2}$
 so that $G(z^2)=g(z).$
\end{lemma}
\proof
Since $g$ is even, it is actually 
holomorphic 
 on $[-R,R]_r$. 
Denoting by $D_r(0):=\{ |z|<r\},$
we have that, since
$g$ is holomorphic and even on $D_r(0),$
 $g(z)=\sum_{j\geq 0} a_{2j} z^{2j},$
 where the power series has a radius of convergence
 $\geq r.$ Then $G(v):=\sum_{j\geq 0} a_{2j}v^j$
 has radius of convergence $\geq r^2.$ 
 It remains to define $G$ in the set
 $\Omega:=[0,R^2]_{r^2}\setminus D_{r^2}(0).$
 Notice  that $\Omega\subset \C\bks(-\infty,0]$; thus, 
 we can define $G(v):=g(\sqrt v)$ for $v	\in \Omega,$
 noting that $z:=\sqrt v\in [0,R]_r.$
Indeed,  if $v\in D_{r^2}(v_0^2),$
 with $v_0\in\R,$ $v_0>r,$ then $\sqrt v\in
 D_r (v_0)$, and this
is equivalent to\footnote{This inclusion follows 
 noting that, for every  $\theta,$
 we have 
 $
 |(v_0 + re^{\ii\theta})^2-v_0^2|\geq r^2.
 $
 The last inequality follows noting that
 it is equivalent to
 $|re^{\ii 2\theta}+2 v_0 e^{\ii\theta}|=
 |re^{\ii \theta}+2 v_0|\geq r$,
 that follows from $v_0>r.$
 }
 $D_{r^2}(v_0^2)\subseteq S(D_r (v_0)),$
 where $S(v):=v^2.$
\eproof

\noindent
Since 
$\cin_\sharp$ is  even in $z,$
by Lemma \ref{geronimo}
we can define  the analytic function
\begin{equation}\label{pachelbel}
\cin_\dag(z^2,\hat p,\sa):=
\cin_\sharp(z,\hat p,\sa)
\end{equation}
which,  by\eqref{rodimento}, satisafies
\begin{equation}\label{3holesbis}
\sup_{v\in (0,\Ro^2)_{\ro^2/16}}  {\modulo}\cin_\dag(v,\hat p,\sa){\modulo}_{\hat D,\ro,\so}\leq \lalla\,.
\end{equation}
Moreover, since 
$v\partial_v \cin_\dag (v,\hat p,\sa)=\frac12 \sqrt v\partial_z
\cin_\sharp (\sqrt v,\hat p,\sa),$
by \eqref{prurito} we get
\begin{equation}\label{3holester}
\sup_{v\in (0,\Ro^2)_{\ro^2/64}}  {\modulo}v\partial_v \cin_\dag
 (v,\cdot,\cdot){\modulo}_{\hat D,\ro,\so}
\leq 8\kappa\lalla\,.
\end{equation}
In the following we will often omit to write the dependence upon
$\hat p$.

\nl
In view of \eqref{palumbo}, \eqref{ummagamma}, 
\eqref{meauxter}, \eqref{4holes} and
\eqref{pachelbel},
we can write
$I(E)=\act^{(2j+1)}_1 (E,\hat p)$ in \equ{4-1}
as 
\begin{eqnarray}
I(E)
&=&
\frac{1}{2\pi}
\int_{\Sa_{{}_{\!\star}}(E )}^{\Sa^{{}^{\!\star}}(E )}
\Big[
\mathcal P
\Big(\sqrt{E-\Gm(\sa)},\sa \Big) 
-
\mathcal P
\Big(-\sqrt{E-\Gm(\sa)},\sa \Big)
\Big]
\, d\sa
\nonumber
\\
&=&
\frac{1}{\pi}
\int_{\Sa_{{}_{\!\star}}(E )}^{\Sa_{{}_{\!\star}}(E )}
\sqrt{E-\Gm(\sa)}\Big(
1
+ \cin_\sharp\big(\sqrt{E-\Gm(\sa)},\sa \big)
\Big)\, d\sa
\nonumber
\\
&=&
\frac{1}{\pi}
\int_{\Sa_{{}_{\!\star}}(E )}^{\Sa^{{}^{\!\star}}(E )}
\sqrt{E-\Gm(\sa)}\Big(
1
+ \cin_\dag\big(E-\Gm(\sa),\sa \big)
\Big)\, d\sa\,.
\label{musicalbox}
\end{eqnarray}
Recalling the definition of 
$\cin_\dag$ in \eqref{pachelbel}, we set
\begin{equation}\label{palidoro}
\tilde \cin(v)=\tilde \cin(v,\hat p,\sa)
:=\cin_\dag(v)+2 v \partial_v \cin_\dag(v)\,,
\end{equation}
which, by \eqref{3holesbis} and \eqref{3holester}, satisfies
\begin{equation}\label{3holesquater}
\sup_{v\in (0,R_0^2)_{r_0^2/64}}  {\modulo}
\tilde \cin (v){\modulo}_{\hat D,r_0,s_0}
\leq 17\kappa\lalla\,.
\end{equation}
 From \eqref{musicalbox} and \eqref{palidoro}
there follows
\beq{lamento2}
\partial_E \act(E)
=
\frac{1}{2\pi}
\int_{\Sa_{{}_{\!\star}}(E )}^{\Sa^{{}^{\!\star}}(E )}
\frac1{\sqrt{E-\Gm(\sa)}}\Big(
1
+ \tilde \cin\big(E-\Gm(\sa),\sa \big)
\Big)\, d\sa\,.
\eeq
Now, note that by \eqref{cimabue} and \eqref{alce} for real values of $\sa$ (and $\hat p$)
$$
\Gm(\sa_{2j+1}+\sa)-\Gm(\sa_{2j+1})
=
\Gm(\sa_{2j+1}+\sa)- E_{2j+1}
\lessdot \suca\sa^2\,.
$$
 Thus, for $E_{2j+1}<E<E_{2j}$ we get
\begin{equation*}\label{lothlorien+}
\frac{1}{\sqrt\suca}\sqrt{E- E_{2j+1}}
\lessdot
\Sa^{{}^{\!\star}}(E)-\sa_{2j+1}\,,\ 
\sa_{2j+1}-\Sa_{{}_{\!\star}}(E)\,.
\end{equation*}
Finally,  by \eqref{lamento2}, \eqref{3holesquater},
we see that
$$
\frac{1}{\sqrt\suca}
\lessdot
\frac{1}{4\pi}
\int_{\Sa_{{}_{\!\star}}(E )}^{\Sa^{{}^{\!\star}}(E )}
\frac1{\sqrt{E-E_{2j+1}}}\, d\sa
\leq
\partial_E \act(E)\,,
$$
proving \eqref{vana}. 
\\
The proof of \eqref{moldavater} is completely  analogous.

\nl
{\bf Step 2:} {\sl  Normal forms close to hyperbolic and elliptic equilibria}

\nl
By Definition \ref{morso}, \equ{alce}, \equ{october}, 
\equ{caviale2}, and Cauchy estimates
one has
\beq{beatoangelico}\ts
\sup_{\hat{p}\in \hat D_\ro}|\Im\sa_{2j}(\hat{p})| \leq
\frac{2\suca\lalla}{\morse \so}
\leq \frac{\so}8\,,\qquad
\sup_{p\in D_{3\ro/4}}|p_1\partial_{\hat{p}}\sa_{2j}(\hat{p})| \leq
(\Ro+\ro)\frac{8\suca\lalla}{\morse \so\ro}
\leq \frac{\so}8\,.
\eeq
Then, the following functions, by \equ{cimabue} and \equ{beatoangelico},
 \beq{caspiterina}
  \cin_*(p, q_1):=\cin(p, q_1+\sa_{2j}(\hat p))\,,\quad
 \Gm_*(\hat p, q_1):=\Gm(\hat p, q_1+\sa_{2j}(\hat p))\,,\quad
  \GO_*( q_1):=\GO( q_1+\bar\sa_{2j})
 \eeq
 satisfy
 \beqa{cimabue2}
 &&\Gm_*(\hat p, 0)=E_{2j}(\hat p)\,,\qquad
 \partial_{q_1}\Gm_*(\hat p, 0)=0\,,
 \nonumber
 \\
  &&\sup_{\hat D_\ro\times \torus^1_{7\so/8}}|\cin_*| \leq
\lalla\,,
\qquad
 \sup_{\torus^1_\so}|\GO_*|\le \suca\,,
\qquad
 \sup_{\hat D_\ro\times \torus^1_{7\so/8}}|\Gm_*-\GO_*| \leq
17\upkappa^3\suca
 \lalla
\,.
\eeqa
 In particular the last estimate follows since for 
 $(\hat p, q_1)\in \hat D_\ro\times \torus^1_{7\so/8}$
 one has
 \beqano\ts
 |\Gm_*(\hat p, q_1)-\GO_*(q_1)|
 &\leq&  
 |\Gm_*(\hat p, q_1)-\GO(q_1+\sa_{2j}(\hat p))|
 +|\GO(q_1+\sa_{2j}(\hat p))-\GO( q_1+\bar\sa_{2j})|
 \\
 &\leq&
\ts \suca
 \lalla
+\frac{8\suca}{\so}
\frac{2\suca\lalla}{\morse \so}
\leq 17\upkappa^3\suca
 \lalla
 \eeqano
 by \equ{october}, \equ{alce} and Cauchy estimates.
Again, by Cauchy estimates, 
 \equ{alce} and \equ{caviale2}
 we get
 \beq{Tarkovski}
  \sup_{\hat p\in \hat D_\ro}|\partial^2_{q_1}\Gm_*(\hat p, 0)-\partial^2_{q_1}\GO_*(0)|
  \leq
  2^6 \upkappa^3\suca\lalla\so^{-2}
  \leq 2^6 \upkappa^6\lalla\morse
  \leq 2^{-10} \morse\,.
 \eeq
By  \equ{A2bis}, 
 $\GO_*$ is
 $\morse$--Morse (and  $\langle \GO_*\rangle=0$);
 in particular it has a maximum at $q_1=0$
 and,  by \equ{ladispoli},  $-\partial^2_{q_1}\GO_*(0)\geq \morse$.
 Recalling \equ{cimabue} and \equ{Tarkovski},  for $\hat p\in\hat D_\ro$, we see that 
\begin{eqnarray}
 &&
\sqrt{\morse/2}\leq
\bar\xxx:=
 \sqrt{-\partial^2_{q_1}\bar\Gm_*(0)/2}
 \leq 
 \sqrt\suca/\so
\,,\qquad
\xxx(\hat p):=
 \sqrt{-\partial^2_{q_1}\Gm_*(\hat p, 0)/2}\,,
 \nonumber
\\
&&
|\xxx(\hat p)-\bar\xxx|
\leq  2^6 \upkappa^6\lalla\sqrt\morse
\leq 2^{-10}\sqrt\morse
\,,\qquad
\frac23\sqrt\morse\leq
|\xxx(\hat p)|\leq 2\upkappa\sqrt\suca\,.
\label{pinturicchio}
\end{eqnarray}
Furthermore,  by \equ{cimabue2}, \equ{caviale2}, 
and by  \equ{pinturicchio}, \equ{cimabue2} and \equ{alce}, we get
for all $\hat p\in\hat D_\ro$:
\beqa{deltoide}\ts
 &&\!\!\!\!\!\!\!\!\!\!\!\!\!\ts \bar\d:=\sqrt{\bar\xxx}\,,\ 
 \d(\hat p):=\frac{\sqrt{\xxx(\hat p)}}{\sqrt[4]{1+\cin_*(0,\hat p,0)}}\,, 
 \qquad
 \frac12\morse^{1/4}\leq
 | \d(\hat p)|
 \leq\upkappa\suca^{1/4}\,,
 \ 
 |\d-\bar\d|\lessdot \lalla\,,\nonumber\\
&&\!\!\!\!\!\!\!\!\!\!\!\!\!\ts
 \bar g:=\bar \xxx
 \,,
  \ 
  g(\hat p):=\sqrt{1+\cin_*(0,\hat p,0)}\xxx(\hat p)\,,
  \ 
  \frac{\sqrt\morse}3\leq
|g(\hat p)|\leq 4\upkappa\sqrt\suca
 \ 
 |g-\bar g |\lessdot \sqrt\suca \lalla\,.
\eeqa
The normal form close to hyperbolic equilibria $(p_1,q_1)=(0,\sa_{2j}(\hat p))$ is detailed in the following  
 
\begin{proposition}\label{bruegel}
There exist positive
constants $\bfco, \bfcu, \bfcd, \bfct$,  depending only on 
 $\upkappa,n$ and 
 satisfying 
 $0<\bfcu<\bfco/8n\bfcd$,
 such that the following holds.
There exist
 a (close to the identity) 
 real analytic symplectic transformation
 \begin{equation}\label{lavedova}
 \Phi_{\rm hp}:(y,x)\in  \{|y_1|<\bfcu \suca^{1/4}\}\times\hat D_{\ro/2}
 \times \{|x_1|<\bfcu \suca^{1/4}\}\times \T^{n-1}_{\so/2} 
 \ \longrightarrow\ 
 (p,q)\in
 \  D_{\ro,\so}
  \end{equation}
and  a function $R_{\rm hp}(z,\hat y)$ with
 \begin{equation}\label{molecolare}
 \sup_{|z|\leq 2\bfcu^2,\ \hat y\in \hat D_{\ro/2}}
|R_{\rm hp}(z,\hat y)|\leq \bfcd\,,\qquad
R_{\rm hp}(0,\hat y)=0\,,\quad
\partial_z R_{\rm hp}(0,\hat y)=0\,,
 \end{equation}
 such that
 \begin{equation}\label{abbacchio}\ts
 \Hpend_{\rm hp}(y,x_1):=\Hpend\circ \Phi_{\rm hp}(y,x)=
 E_{2j}(\hat y)+g(\hat y) (y_1^2-x_1^2)
+\suca R_{\rm hp}\left(\frac{y_1^2-x_1^2}{\sqrt\suca} , \hat y\right)\,.
\end{equation}
Moreover $\Phi_{\rm hp}$ has the form
 \beqa{pesante}
 &&p_1= 
 \d(\hat y)\Big(y_1+\suca^{1/4} a_1(\suca^{-1/4}y_1,\hat y,\suca^{-1/4}x_1)\Big)\,,\quad 
 \hat p=\hat y,
  \\
 &&q_1=
 \sa_{2j}(\hat y)+\frac{1}{\d(\hat y)}
\Big( x_1+
 \suca^{1/4} a_2(\suca^{-1/4}y_1,\hat y,\suca^{-1/4}x_1)\Big)\,,
 \quad 
 \hat q=\hat x+\hat a(y,x_1)\,,
 \nonumber
 \eeqa
for suitable holomorphic functions
$a_1,a_2,\hat a$,
such that
\beq{chedolore}
	\sup_{W_{\!\bfco,\ro}}
		 |a_i|\le \bfcd\,,\qquad 
W_{\!\bfco,\ro}:=\{\ts|\tilde y_1|<{\bfco}/2 \}\times\hat D_{\ro/2}  \times \{\ts|\tilde x_1|<{\bfco}/2 \}\,,
	\eeq
	and are
	at least quadratic in $\tilde y_1,\tilde x_1$.
	Moreover denoting $\bar R_{\rm hp}:=R_{\rm hp}|_{\lalla=0}$ we have
	\begin{equation}\label{biretta}
|R_{\rm hp}-	\bar R_{\rm hp}|=O(\lalla)\,. 
\end{equation}
Finally,  for every $\hat y \in \hat D_{\ro/2}$,
the image of  the restriction  of the map in \equ{pesante}
$$
(y_1,x_1)\in \{|y_1|<\bfcu \suca^{1/4}\}
 \times \{|x_1|<\bfcu \suca^{1/4}\}
 \to 
 (p_1,q_1)
$$
 contains the (complex) set
 \beq{brodino}
\{|p_1|\leq2\bfct \sqrt\suca\}
 \times \{|q_1-
 \sa_{2j}(\hat p)|\leq2\bfct \}\,.
 \eeq
\end{proposition}
The normal form close to elliptic equilibrium $(p_1,q_1)=(0,\sa_{2j+1}(\hat p))$ is detailed in the following  

 \begin{proposition}\label{bruegelell}
There exist
 a (close to the identity)
  real analytic 
  symplectic transformation
 $\Phi_{\rm el}$ as in \eqref{lavedova} and \eqref{pesante}
 and  a function $R_{\rm el}(z,\hat y)$ 
 as in \eqref{molecolare} and \eqref{biretta}
  such that
 \begin{equation}\label{abbacchioell}\ts
 \Hpend_{\rm el}(y,x_1):=\Hpend\circ \Phi_{\rm el}(y,x)=
 E_{2j+1}(\hat y)+g(\hat y) (y_1^2+x_1^2)
+\suca R_{\rm el}\left(\frac{y_1^2+x_1^2}{\sqrt\suca} , \hat y\right)\,.
\end{equation}
Finally,  for every $\hat y \in \hat D_{\ro/2}$,
the image of  the restriction 
$$
(y_1,x_1)\in 
\{|y_1|<\bfcu \suca^{1/4}\}
 \times \{|x_1|<\bfcu \suca^{1/4}\}
 \to 
 (p_1,q_1)
$$
 contains the (complex) set
 \beq{brodinoell}
\{|p_1|\leq2\bfct \sqrt\suca\}
 \times \{|q_1-
 \sa_{2j+1}(\hat p)|\leq2\bfct \}\,.
 \eeq
\end{proposition}
The   proof of Proposition~\ref{bruegel} is rather standard;  
for completeness it is included in Appendix~\ref{oceania}. \\
The proof of Proposition~\ref{bruegelell} is completely analogous\footnote{Indeed, in one dimension,  from a complex point of view,
the Birkhoff normal form is the same both in the
hyperbolic and in the elliptic case.
} and is omitted.


\nl
{\bf Step 3:} {\sl  The action  functions close to the elliptic equilibrium}

\nl
Setting $I_1:=(y_1^2+x_1^2)/2$ and $\hat I:=\hat y$,
we have that, by Proposition 
\ref{bruegelell}, the function 
\beq{abbacchioell2}\ts
{\mathtt E}(I):= E_{2j+1}(\hat I)+2g(\hat I) I_1
+\suca R_{\rm el}\left(\frac{2I_1}{\sqrt\suca} , \hat I\right)
\eeq
is well defined and holomorphic for
$|I_1|\leq \bfcu^2\sqrt\suca/2$ and 
$\hat I \in \hat D_{\ro/2}$.
Since $R_{\rm el}$ is (at least) quadratic
in its first entry, by \eqref{deltoide} and
\eqref{alce}, recalling
\eqref{mortazza}, \eqref{october}
and taking $\cc$ large enough in \eqref{caviale2}, we see that, for   a 
suitable constant $0<\bfcq<\bfct$ 
 depending only on $n$ and $\upkappa$,
we can invert the expression
$\mathtt E(I)=E$ finding $I_1(E,\hat I)$
 which solves
 $$
 \mathtt E\big(I_1(E,\hat I),\hat I\big)=E\,,
 \qquad \mbox{for}\ \ \ 
 |E-\bar E_-|<\bfcq \suca\,,\ \ \ 
 \hat I \in \hat D_{\ro/2}\,.
 $$
 It turns out that the function $I_1$ above
is exactly the action function
introduced in \eqref{4-1}.
Indeed, since the map $ \Phi_{\rm el}$
in \eqref{abbacchioell}
is symplectic,
 for every $ |E-\bar E_-|<\bfcq \suca$ the area enclosed by the level curve
$\g_{2j+1}(E;\hat I)$ in \eqref{poirot} is equal to the one
included
by the level curve $\Hpend_{\rm el} =E$, which
is simply the circle $\frac12(x_1^2+y_1^2)=I_1(E).$

\nl
Hence, 
formula \eqref{lamponek} (i.e., the analyticity of action as a function
of energy close to a minimum
and $I_1(E_-)=0$)
 follows.


\nl
{\bf Step 4:} {\sl  Away from the elliptic equilibrium}

\nl
Here, 
we will often omit to write the dependence on 
the dumb actions $\hat p=\hat I=\hat y$.

\nl
Let us consider  the action   function $\act(E)$ 
 defined in \eqref{musicalbox} for
\beq{ano}
 |E-\bar E_-|\geq\bfcq \suca\,.
 \eeq
By Green's Theorem,
\begin{equation}\label{ziogianni}
\act(E)=\frac1{2\pi}\int_{\Omega(E)}
dq_1 dp_1\,,
\end{equation}
where $\Omega(E)$ is the bounded  portion of plane encircled by the curve
$\g_{2j+1}(E)$  defined in \equ{poirot}, namely\footnote{The definition of $\mathcal P$  is given in Lemma~\ref{amatriciana}.}
\begin{eqnarray}\label{ciocco}
&&\Omega(E)=
\Big\{(p_1,q_1)\ |\ \Sa_{{}_{\!\star}}(E )<q_1<\Sa^{{}^{\!\star}}(E )\,,\ \ 
\\
&&\qquad\qquad
\ts \mathcal P
\Big(-\sqrt{E-\Gm(q_1)},q_1 \Big) 
\leq
p_1
\leq
\mathcal P
\Big(\sqrt{E-\Gm(q_1)},q_1 \Big)
\Big\}\,.
\nonumber
\end{eqnarray}
Consider first the case in which
 we are away also from the hyperbolic equilibrium.
By \eqref{musicalbox}
we have
$$
\int_{\Omega(E)} \ dq_1dp_1=2
\int_{\Sa_{{}_{\!\star}}(E )}^{\Sa^{{}^{\!\star}}(E )}
\sqrt{E-\Gm(\sa)}\Big(
1
+ \cin_\dag\big(E-\Gm(\sa),\sa \big)
\Big)\, d\sa\,,
$$
which contributes to $I(E)$
with a holomorphic\footnote{Notice that there is no problem 
in  $\Sa_{{}_{\!\star}}(E )$, $\Sa^{{}^{\!\star}}(E )$
where the square root vanishes. 
Actually,  close to these points it is convenient to 
write  $\Omega(E)$ as a normal set 
with respect to $q_1$ and not to $p_1$.} and bounded (by some
constant depending only on $\upkappa$
and $n$) term.

\nl
Let us finally consider the case
close to the hyperbolic point.
Recalling \eqref{abbacchio}
let us consider the equation
\begin{equation}\label{maia}\ts
E_{2j}(\hat y)+
g(\hat y) J
+\suca R_{\rm hp}\left(\frac{J}{\sqrt\suca} , \hat y\right)=E
\end{equation}
By the inverse function theorem we construct
a holomorphic function $F( z,\hat y)$
with 
\begin{equation}\label{carciofino}
\sup_{| z|<\bfcc, \hat D_{\ro/2}}
|F(z, \hat y)|\leq 1/2\bfcc\,,
\end{equation}
for $0<\bfcc<\bfcq<\bfct$ small enough depending only on $\upkappa,n$,  such that
the equation in \eqref{maia} is solved by\footnote{For
real values of $\hat y$ and $E$ we are in the case $E<E_{2j}(\hat y)$,
namely $z>0$.}
\begin{equation}\label{ziatita}\ts
J(E,\hat y):=-\sqrt\suca
  J_{\rm hp}\left(\frac{E_{2j}(\hat y)-E}{\suca} ,\hat y\right)\,,\qquad
  J_{\rm hp}(z,\hat y):=
\frac{\sqrt\suca }{g(\hat y)}z\left(1+
  z F( z,\hat y)\right)\,,
\end{equation}
where $  J_{\rm hp}$ solves the equation 
\begin{equation}\label{diamante}
\frac{g(\hat y)}{\sqrt\suca}   J_{\rm hp}-R_{\rm hp}(-  J_{\rm hp},\hat y)= z\,.
\end{equation}
\nl
Recalling \eqref{ziatita} we set
\begin{equation}\label{satizza}
x_1(E,y):=\sqrt{-J(E,\hat y)+y_1^2}
\end{equation}
then, by \eqref{maia} and \eqref{abbacchio},
\begin{equation}\label{grandinata}
 \Hpend_{\rm hp}(y,x_1(E,y))\equiv E\,.
\end{equation}
Fix
$
\bar x_1:=\bfcc \suca^{1/4}
$
and  define
\beq{herrmannelig}
\bar y_1=\bar y_1(E,\hat y):=\sqrt{J(E,\hat y)+\bar x_1^2}\,.
\eeq
Consider the holomorphic functions 
$q_1=q_1(y,x_1)$ and $p_1=p_1(y,x_1)$
defined in  \eqref{pesante}.
Set
\beq{bach}
\bar p_1^\pm(E,\hat y):=
p_1\big(\pm\bar y_1(E,\hat y),\hat y,\bar x_1\big)\,,
\qquad
\bar q_1^\pm(E,\hat y):=
q_1\big(\pm\bar y_1(E,\hat y),\hat y,\bar x_1\big)\,.
\eeq
Note that\footnote{Omitting $\hat y$.}
\beq{panico}
\bar p_1^\pm(E)=
\mathcal P
\ts\Big(\pm\sqrt{E-\Gm(\bar q_1^\pm(E))},\bar q_1^\pm(E) \Big)
\eeq
For every fixed $\hat y$ (that we will omit to write)
we invert the expression
$
p_1=p_1(y_1,\bar x_1)
$, with $|y|<\bfcu\suca^{1/4}$,
finding a holomorphic function $\tilde y_1(p_1)$
such that
\beq{caifa}
p_1=p_1\big(\tilde y_1(p_1),\bar x_1\big)\,,\qquad
\bar y_1^\pm(E)=\tilde y_1(\bar p_1^\pm(E))
\,.
\eeq
Set $\tilde q_1(p_1):=q_1(\tilde y_1(p_1),\bar x_1)$.
For real value of $E$ and $\hat p$, we split the 
integral in \eqref{ziogianni} in two parts
\beq{neve0}
\int_{\Omega(E)} \ dq_1dp_1
=
\int_{\Omega_1(E)} \ dq_1dp_1
+\int_{\Omega_2(E)} \ dq_1dp_1
\eeq
where
$$
\Omega_1(E) :=
\{(p_1,q_1)\in 
\Omega(E)\ |\ q_1\leq \tilde q_1(p_1)\,,
\ \bar p_1^-\leq p_1\leq \bar p_1^+\}
\,,
\quad
\Omega_2(E)
:=\Omega(E)\setminus\Omega_2(E)
\,.
$$
As above the term
$
\int_{\Omega_2(E)} \ dq_1dp_1$
, contributes to $I(E)$
with a holomorphic and bounded term.

\nl 
Recalling \eqref{abbacchio} and 
\eqref{brodino} and setting
$\tilde \Omega_1(E):=\Phi_{\rm hp}^{-1}\big(
\Omega_1(E)\big)$ we have that
\beq{neve}
\int_{\Omega_1(E)} \ dq_1dp_1
=\int_{\tilde\Omega_1(E)} \ dx_1dy_1\,.
\eeq
Note that by the above construction
$$
\tilde\Omega_1(E)=\{
(y_1,x_1)\ :\ -\bar y_1(E)\leq y_1\leq \bar y_1(E)\,,
\ 
x_1(E,y_1)\leq x_1\leq \bar x_1
\}\,,
$$
then
\beq{neve2}
\int_{\tilde\Omega_1(E)} \ dx_1dy_1
=2\bar x_1 \bar y_1(E)-\int_{-\bar y_1}^{\bar y_1}
 x_1(E,y_1)\, dy_1\,.
\eeq
 On the other hand by \eqref{satizza} and \eqref{herrmannelig}
\begin{eqnarray*} 
 &&\int_{-\bar y_1}^{\bar y_1}
 x_1(E,y_1)\, dy_1
=2 \int_0^{\bar y_1}
\textstyle
\sqrt{-J(E)+y_1^2}\, dy_1
\\
&&
\textstyle
=\bar y_1\sqrt{-J(E)+\bar y_1^2}
-J(E)\Big(
\ln\frac{\sqrt{-J(E)+\bar y_1^2}+\bar y_1}{\suca^{1/4}}
-
\ln\frac{\sqrt{-J(E)}}{\suca^{1/4}}
\Big)
\\
&&
\ts
=\bar y_1 \bar x_1
\ -J(E)
\ln\frac{\bar x_1+\bar y_1}{\suca^{1/4}}
\, +\frac12 J(E)\ln\frac{-J(E)}{\sqrt\suca} \,.
 \end{eqnarray*}
 Note that the first two terms above are 
 holomorphic functions of 
 $E$ up to $E=0$,
instead the term 
 $\textstyle\frac12 J(E)\ln (-J(E)/\sqrt\suca)  $
contains the singular term.
Recalling \eqref{mortazza} and setting
$z=(E_{2j}-E)/\suca=(E_+ -E)/\suca$, 
 the last term is transformed into (recalling \eqref{ziatita})
 \begin{eqnarray*}\ts
  \frac12 J(E)\ln\frac{-J(E)}{\sqrt\suca} 
  &=&\ts
  -\frac{\sqrt\suca} {2}  J_{\rm hp}(z)\ln  J_{\rm hp}(z)
  \\
 &=&\ts
  -\frac{\sqrt\suca} {2}  J_{\rm hp}(z)\ln\left( 
  \frac{\sqrt\suca }{g(\hat y)}\left(1+
   z F( z,\hat y)\right)\right)
  -\frac{\sqrt\suca} {2}  J_{\rm hp}(z)\ln z\,,
 \end{eqnarray*}
 where the last term is the singular one, namely
 $$
 \ts -\frac{\sqrt\suca} {2}  J_{\rm hp}(z)\ln z
 =-\frac{\suca }{2 g(\hat y)}\left(1+
 z  F( z,\hat y)\right)\, z \ln z\,.
 $$
 Recalling \eqref{ziogianni}, \eqref{neve0}, \eqref{neve} and \eqref{neve2}, this implies that 
 the singular term in $I(E_+ -\suca z)$ in
  \eqref{LEGOk} is
  $$
  \ts
  \frac{\suca }{4\pi g(\hat y)}\left(1+
  z F( z,\hat y)\right)\, z \ln z\,,
  $$
namely
 \begin{equation}\label{ceci}\ts
 \psi^{(2j+1)}_+(z)
 =
 \frac{\suca }{4\pi g(\hat y)}\left(1+
 z  F( z,\hat y)\right)\,.
 \end{equation}
 This proves \eqref{LEGOk}. 
 
 \nl
By taking $\cc$ in \eqref{caviale2} large enough, by \eqref{carciofino} and \eqref{deltoide}
 the first estimate in \eqref{pappagallok} and  
 the first and third estimates in
 \eqref{ciofecak} follow for
 $\psi^{(2j+1)}_+$.

\nl 
 Now consider the corresponding functions when $\lalla=0$ (namely, $\bar I, \bar J_{\rm hp}, \bar \psi^{(2j+1)}_+$, etc.).
 Observe, in particular\footnote{Recall \eqref{ziatita}.}, that 
 $$\ts 
 \bar J_{\rm hp}(z):=
\frac{\sqrt\suca }{\bar g}z\left(1+
  z \bar F(z)\right)\,,
  $$
 is the solution of the equation
 \begin{equation}\label{diamante2}\ts
\frac{1}{\sqrt\suca} \, \bar g\bar J_{\rm hp}-\bar R_{\rm hp}(-\bar J_{\rm hp})= z\,.
\end{equation}
 corresponding to \eqref{diamante}.
 Then recalling \eqref{ziatita}
 \beqano\ts
 z\left(1+
  z \bar F(z)\right)
  -\bar R_{\rm hp}(-\bar J_{\rm hp})&=&\ts
 \frac{\bar g}{\sqrt\suca} \bar J_{\rm hp}-\bar R_{\rm hp}(-\bar J_{\rm hp})
 =
 \frac{g(\hat y)}{\sqrt\suca}   J_{\rm hp}-R_{\rm hp}(-  J_{\rm hp},\hat y)
\\
&=&
 z\left(1+
  z  F(z,\hat y)\right)-R_{\rm hp}(-  J_{\rm hp},\hat y)\,.
 \eeqano
 By \eqref{deltoide} and \eqref{biretta}, we get
 \begin{equation}\label{pippero}
 |F-\bar F|\lessdot \lalla\,,\qquad
|  J_{\rm hp}-\bar J_{\rm hp}|\lessdot \lalla\,.
\end{equation}
Since the unperturbed singular term is
$$\ts
 \bar\psi^{(2j+1)}_+(z)= 
 \frac{\suca }{4\pi \bar g}\left(1+
  z\bar F( z)\right)
 $$
 by \eqref{ceci}, \eqref{deltoide} and \eqref{pippero}
 we get the second estimate in \equ{trippa} in the $+$ case when
 $i=2j+1$; 
 the other cases are analogous.
 \\
Since $\bar\psi^{(2j+1)}_+$ is independent of $\hat\act$,
by \eqref{trippa} and Cauchy estimates, we get
  the second estimate in \eqref{pappagallok} for
 $\psi^{(2j+1)}_+$; the other estimates are analogous.  
 
 \nl
 It remains to prove \eqref{rosettaTH}.
 In proving \eqref{rosettaTH} we consider only the 
 crucial zone close to maximal energies;
 in particular we consider the domain
\beq{ambe}
{\mathcal E}^{i}_\loge\cap\{|E-\bar E_+|\leq \suca/2\cc\}\,.
\eeq
Indeed 
 in the other parts the estimates are simpler 
 and can be directly derived from the representation
 formula \eqref{frascati2}, \eqref{lamento2}
 and the estimate \eqref{3holesquater};
 noting also that, by \eqref{ano} and \equ{cimabue}, 
 the function $\Sa_{{}_{\!\star}}(E)$, resp. $\Sa^{{}^{\!\star}}(E)$, and 
 $\bar\Sa_{{}_{\!\star}}(E)$,  resp. $\bar\Sa^{{}^{\!\star}}(E)$, 
 are close:
 $$
 |\Sa_{{}_{\!\star}}(E)-\bar\Sa_{{}_{\!\star}}(E)|\,,\ 
 |\Sa^{{}^{\!\star}}(E)-\bar\Sa^{{}^{\!\star}}(E)|\, \lessdot\, 
 \suca\lalla\,.
 $$
Let us consider the domain in \eqref{ambe}, where
we can use the representation \eqref{LEGOk} and estimates in
 \eqref{pappagallok}--\eqref{trippa}.
 The first and second estimate in \eqref{rosettaTH} 
 directly follow from \eqref{pappagallok}.
 Let us now consider the third estimate in \eqref{rosettaTH}.
 Denote
\beqno
\left\{
\begin{array}{l}
z:=(E_{2j}-E)/\suca \,,\quad 
 z_*:=(\bar E_{2j}-E)/\suca \,,\\
f(z):=\phi'(z)+\psi'(z)z\ln z+\psi(z)(1+\ln z)\,,\\
\bar f(z):=\bar \phi'(z)+\bar\psi'(z)z\ln z+\bar\psi(z)(1+\ln z)
\end{array}\right.
 \eeqno
 and observe that, by \eqref{october},  $|z-z_*|=|E_{2j}-\bar E_{2j}|/\b\lessdot \lalla$. Then, recalling \eqref{ambe},
$$
\suca\big|
\partial_E \act_1^{(2j)}(E)
-
\partial_E \bar \act_1^{(2j)}(E)\big|=|f(z)-\bar f(z_*)|\lessdot
{\sqrt\suca\lalla}/{\loge}
$$
 by \eqref{pappagallok}, \eqref{trippa} and Cauchy estimates.
 The proof of Theorem~\ref{glicemiak} is complete.
 \qed 


\section{The complex Arnol'd--Liouville transformation}

In this  section we discuss the complex properties (including analyticity radii) of the Arnol'd--Liouville transformation, which allow, in particular,  to give upper bound on the derivatives of the energy functions in complex domain.

\nl
For every fixed $\hat p\in\hat D$,
given $\act_1^{(i)}(E,\hat p)$ as in \equ{4-1}
the  action function $(p_1,q_1)\to \act_1^{(i)}(\Hpend(p,q_1),\hat p)$
can be symplectically completed\footnote{Uniquely fixing, e.g., 
$\f^{(i)}_1(p,0)=0$.} with the angular term
$(p_1,q_1)\to \f^{(i)}_1(p,q_1)$.
We shall call $\cFiq=\cFiq(\act,\ang_1)$
the inverse of the map 
$$
(p,q_1)\to (\act,\f_1):=
\big(\act_1^{(i)}(\Hpend(p,q_1),\hat p),\hat p,\f^{(i)}_1(p,q_1)\big)\,.
$$
Note that 
the Arnol'd-Liouville 
`suspended'
symplectic transformation
$\cFiq$ integrates  $\Hpend$, i.e.,
\begin{equation}\label{red2}
\Hpend\circ\cFiq (\act,\ang_1)
={\mathtt E}^{(i)} (\act)\,,
\end{equation}
where ${\mathtt E}^{(i)}$
is the inverse of 
$\act_1^{(i)}$, namely
\beq{ibrido}
{\mathtt E}^{(i)}
\big(\act_1^{(i)}(E,\hat \act),\hat \act\big)=E\,.
\eeq

\nl
Next, we introduce suitable decreasing subdomains 
$\Bu^i (\loge)$ of $\Bu^i$ 
depending on a non negative parameter $\loge$ so that $\Bu^i(0)=\Bu^i$  and 
such that 
the map $\Fiq$ 
has, for positive $\loge$,  a holomorphic extension
on a suitable complex neighborhood
of $\Bu^i(\loge)\times\T^n$.

\nl
Define
\begin{equation}\label{sinistro}
\logemax=\logemax(\hat \act):=\big(E_+(\hat \act)-E_-(\hat \act)\big)/\suca\,,\qquad
\blogemax:=
\big(\bar E_+-\bar E_-\big)/\suca
\,.
\end{equation} 
Notice  that, by \equ{alce}, the definition of $\morse$ and \equ{cimabue}  one has 
\beq{latooscuro}
\frac{1}{\upkappa}\leq 
\frac{\morse}{\suca}
\leq
\blogemax
\leq
2\,;
\eeq
notice also that, by  \eqref{october},  we have 
\begin{equation}\label{destro}
|\logemax-\blogemax|\leq 6\upkappa^3\lalla\,,
\end{equation}
so that, since $\lalla\leq 1/\cc^2$ and  $\cc \geq  2^8\upkappa^3$ (compare Theorem~\ref{glicemiak}), one has 
\begin{equation}\label{sinistrorso}
\logemax\geq 1/2\upkappa\,.
\end{equation} 
Next, for   $0\le\loge\leq \logemax$
define:
\begin{eqnarray}\label{arista}
&&\acci^{(2j-1)}_- (\hat\act;\loge):=0\,,\ \forall 1\leq j\leq N\,,\
\nonumber\\
&& \acci^{(2j)}_-(\hat\act;\loge):=
\act_1^{(2j)}\big(E^{(2j)}_-(\hat \act)+\loge\suca,\hat \act\big)\,,\ 
\forall 0\leq j\leq N,
\nonumber
\\
&&\acci^{(i)}_+(\hat\act;\loge):=
\act_1^{(i)}\big(E^{(i)}_+(\hat \act)-\loge\suca,\hat \act\big)\,,\ \forall
0< i< 2N\,,
\nonumber
\\
&&\acci^{(i)}_+(\hat\act;\loge):=
\act_1^{(i)}\big(\Ro^2+  \Ro\ro,\hat \act\big)\,,\ i=0,2N\,,
\end{eqnarray}
and, for $0\leq i\leq 2N$,
 \begin{equation}\label{playmobilk}
\Bu^i (\loge):=\Big\{ \act=(\act_1,\hat \act)\ |\ \hat \act\in\hat D,\ \  \acci^{(i)}_-(\hat \act;\loge)<\act_1
<\acci^{(i)}_+(\hat \act;\loge)
\Big\}\subseteq  \R^n\,,\  \  \Bu^i:=\Bu^i (0)\,. 
\end{equation}
Note that\footnote{Recall \eqref{cimabue}--\eqref{alce}.}
\begin{equation}\label{pratola}
 \diam \Bu^i(0)
\leq 2\big(\Ro +\diam\hat D)\,,
\qquad \forall\ \, 0\leq i\leq 2N\,.
\end{equation}
Setting\footnote{Recall
\eqref{insalatina} and \eqref{gonatak}.}
$$
\check\cM^i:=\{ (p,q_1)\in\R^n\times\T\ \ {\rm s.t.}\ \ 
 (p_1,q_1)\in\cM^i(\hat p)\,,\ \hat p\in\hat D \}
$$
we have
\begin{equation}\label{peppapig}
\check\cM^i
=\check\cM^i(0)
= \cFiq (\Bu^i \times\T)
=\bigcup_{0< \loge\leq 1/\cc}\check\cM^i(\loge)\,,
\quad
{\rm where}
\quad
\check\cM^i (\loge)
 :=
 \cFiq (\Bu^i (\loge)\times\T)\,.
\end{equation}

\begin{theorem}\label{barbabarba}
Under the  hypotheses of Theorem~\ref{glicemiak}
there exists 
$\hcc=\hcc(n,\upkappa)\ge 4 \cc^2$ 
depending only on
$n$ and $\upkappa$ such that,
taking
\beq{caviale3}
\lalla\leq 1/\hcc\,,
\eeq
 for any
$0\leq i\leq 2N$, the symplectic transformation $\cFiq$
 extends, for any 
 $0< \loge\leq 1/\hcc$,
to a real--analytic map 
\begin{equation}\label{pediatra2bis}
\Fiq :\big(\Bu^i (\loge)\big)_{\!\rho}\times\T^n_\s\ \to\ 
D_{\ro}\times\T^n_{\so/4}
\end{equation}
 with 
\begin{equation}\label{blueeyes}\ts
\rho	=\frac{\sqrt\suca} {\hcc}\, \loge |\log \loge|  \,,\qquad \s= \frac1{\hcc|\log \loge|}\,.  
\end{equation}
Moreover, on $\big(\Bu^i (\loge)\big)_{\!\rho}$ we have
 \begin{eqnarray}\label{ofena}
 &&\ts
 \big|\partial_{\act_1} \mathtt E^i \big|
\leq 
\hcc\sqrt{\suca+|\mathtt E^i|}\,, 
  \quad
\big|\partial^2_{\act_1} \mathtt E^i \big|
\leq 
\frac{\hcc}{\hat\loge}\,,\qquad\qquad\qquad(\hat\loge:=\loge|\ln\loge|^3)
\nonumber
\\
&&
\ts
\big|\partial^2_{\act_1 \hat \act} \mathtt E^i \big|
\leq
\hcc\frac{\lella}{\hat\loge}
\,,
\quad
\big|\partial^2_{\hat \act} \mathtt E^i \big|
\leq \hcc
\Big(
\frac{\sqrt\suca}{\ro}\act_1^{i}
+
 \frac{\lella}{\hat\loge}
\Big)
\lella\,.
\end{eqnarray}
Finally, 
 we have
  \begin{equation}\label{inthecourtk}
\meas\Big(
\big(D^\flat\times\T\big)\ \setminus\ 
\bigcup_{0\leq i\leq 2N}\check\cM^i (\loge)\Big)
\leq \hcc\, \sqrt\suca
\meas(\hat D)\ \loge |\log \loge|\,,
\end{equation}
 where $D^\flat:=
  (-\Ro -\ro/3,\Ro+\ro/3)\times\hat D.$
 \end{theorem}
 
\rem
\label{icm2018}
(i) The complete symplectic action--angle map 
$\Fiq: (\act,\ang)\to (p,q)$ has the form  
\beq{bruford2} 
\Fiq(\act,\ang)=\casitwo{(\upeta^i,\hat \act, \uppsi^i,\hat \ang+\upchi^i)\,,}
{0<i<2N\,,}
{(\upeta^i,\hat \act, \ang_1+ \uppsi^i,\hat \ang+\upchi^i)\,,}
{i=0,2N\,,}
\eeq
where $\upeta^i ,\upchi^i ,\uppsi^i $ are function of $(\act,\ang_1)$ only and are 
$2\pi$--periodic in $\ang_1,$
and, in the case $i=0,2N,$
$\sup|\partial_{\ang_1}\uppsi^i |<1$. \\
Notice   that, since
$\bfcu\geq  4\upkappa\geq 16$, by
\eqref{alce} and \eqref{blueeyes} we get
$\rho\leq 2^{-8}\ro$ and
 $\s\leq \so/4.$
By \eqref{pediatra2bis} we also get
$|\Im\uppsi^i|_{\rho,\s}\leq \so/2$
for every $0\leq i\leq 2N.$
Analogously\footnote{Actually a better estimate holds: it is smaller
than some constant by 
$\lella\so$, where $\lella$ was defined in \eqref{pappagallok}.}
$|\Im\upchi^i_j|_{\rho,\s}\leq \so/2$ for every $j=2,\ldots,n$.

\nl
(ii)
Notice the different topologies of this map:  For $1\le i\le 2N-1$ the motion is {\sl librational}, i.e., the 
$q_1$--coordinate oscillates around relative (stable) equilibria, while for $i=0$ and $i=2N$ the motion is {\sl rotational}, 
corresponding to the $q_1$--coordinate rotating  in the unbounded regions of phase space `outside' separatrices; such regions correspond to the labels $i=2N$ (upper unbounded region) and  $i=0$ (lower unbounded region). 

\nl
(iii) For related estimates on the analyticity strip in the angles, see \cite{Oku}. 
\erem

\proof {\bf  of Theorem~\ref{barbabarba}}
The fact that the map $\cFiq$
 extends a complete 
 symplectic transformation $\Fiq$
 directly follows by the Arnold--Liouville
 Theorem. Here we have only 
 to evaluate the analyticity radia.
 For brevity we will often drop the index $i$ and the dumb actions as well as
  will often write
 $I(E)$ instead of $I_1^{(i)}(E)$.
 As above we will consider only the case
 $i$ odd, the other one being similar.

\nl
In order to prove \eqref{pediatra2bis}--\eqref{blueeyes} we introduce energy-time $(E,t)$ 
(symplectic) coordinates,
which are a simple rescaling of action-angle variables $(I,\f)$. Indeed
considering the integrable hamiltonian ${\mathtt E} (\act)=E$
we have 
that the action  and the angular velocity are constant $\dot\f=
\partial_\act {\mathtt E} (\act)$ so that $\f(t)=\partial_\act {\mathtt E} (\act) t$
or, using $(E,t)$  as independent variables
\begin{equation}\label{alice}
I=I(E)\,,\qquad
\f=\f(E,t)=\frac{t}{\partial_E \act(E)}\,.
\end{equation}
We restrict to the zone around hyperbolic points where one has worst estimates.
Here we first pass in $(y_1,x_1)$ coordinates obtaining the Hamiltonian 
$\Hpend_{\rm hp}(y_1,x_1)$
in 
\eqref{abbacchio}.
Secondly we pass to coordinates $(E,t)$ setting
\beq{tirzan}
E:=\Hpend_{\rm hp}(y_1,x_1)\,,\qquad
t:=\frac1{w(E)} {\rm arctanh} 
\left(\frac{y_1}{\sqrt{-J(E)+y_1^2}}\right)\,,
\eeq
where
$w(E):=2 g \big(1+ \sqrt\suca \partial_z R_{\rm hp}(J(E)/\sqrt\suca ) \big)$.
Indeed, by the Hamilton equations for $\Hpend_{\rm hp}$  we leads to 
$$
\dot y_1=-\partial_{x_1}\Hpend_{\rm hp}(y_1,x_1)
=2 g x_1\Big(1+ \sqrt\suca \partial_z R_{\rm hp}\big((y_1^2-x_1^2)/{\sqrt\suca} \big) \Big)
=w(E)\sqrt{-J(E)+y_1^2}
\,,
$$
which can be easily
integrated (by separation of variables) giving the expression for the time in
\eqref{tirzan}.

\begin{lemma}
There exists a small constant $0<\bfcs\leq \min\{\bfcu^2,1\}/2$
depending only on $\upkappa$ and $n$
such that, taking
 \beq{ferrari}
0<\loge\leq \bfcs\quad \mbox{and}\quad
\tilde E:=E_+-\suca\loge\,,
 \eeq
 the map 
 $
(y_1,x_1) \in  \{|y_1|<\bfcu \suca^{1/4}\}
 \times \{|x_1|<\bfcu \suca^{1/4}\}
 \ \to \ (E,t)
 $
 in \eqref{tirzan} is invertible for
 \beq{ninetto}
 |E-\tilde E|\leq\bfcs\suca\loge\,,
 \qquad |t|<\bfcs /\sqrt\suca\,.
 \eeq
\end{lemma}
\proof
Inverting the second expression in \eqref{tirzan}
we get
\beq{enna}
y_1(E,t)=\sqrt{-J(E)}\sinh(w(E)t)\,.
\eeq
Moreover \eqref{ziatita},
\eqref{satizza} and \eqref{grandinata}
we have
\beq{enna2}
x_1(E,t)=\sqrt{-J(E)+\big(y_1(E,t)\big)^2}\,.
\eeq
We have to check that the above functions
$y_1$ and $x_1$ are defined on the set in 
\eqref{ninetto}.
Set
\begin{equation}\label{tampone}
 \sqrt\suca \loge\lessdot J_*:=\frac{\suca} {\bar g}\loge\lessdot \sqrt\suca \loge\,,
\end{equation}
with
$\bar g$ defined in \eqref{deltoide}
(recall also \eqref{alce}).
Recalling the definition of $J(E)$ in \eqref{ziatita}
we have 
$$
J(\tilde E)=-\frac{\suca} {g}\loge(1+\loge F(\loge))
$$
and, by \eqref{deltoide} and \eqref{carciofino}
$$
|J(\tilde E)+J_*|
=\frac{\suca} {\bar g}\loge
\left|
\frac{\bar g}{g}(1+\loge F(\loge))-1
\right|
\lessdot \sqrt\suca\loge(\lalla+\loge)\,.
$$
Finally, since by \eqref{ziatita},\eqref{carciofino} and \eqref{deltoide} $|\partial_E J|\lessdot1/ \sqrt\suca$ we have,
for $|E-\tilde E|\leq \bfcs\suca\loge$,
$$
|J(E)-J(\tilde E)|\lessdot \bfcs\sqrt\suca\loge\,,
$$
then
we get
\begin{equation}\label{tampone2}
|J(E)+J_*|\lessdot (\bfcs+\lalla)\sqrt\suca\loge.
\end{equation}
By \eqref{tampone}
and
taking $\bfcs$ small enough and $\hcc$
large enough, 
we obtain
\begin{equation}\label{schiaffino}
-\Re J(E)\geq 4\bfcs \sqrt\suca\loge\,,
\end{equation}
for any $|E-\tilde E|\leq {\bf c}\suca\loge$.
Moreover
\beq{cipster}
|w(E)|\stackrel{\eqref{deltoide},\eqref{molecolare}}\lessdot \sqrt\suca\,,
\qquad
|J(E)|\stackrel{\eqref{tampone2},\eqref{tampone}}\lessdot \sqrt\suca\loge\,.
\eeq
Then, recalling \eqref{enna},
we get
$$
|y_1(E,t)|\leq \sqrt{\bfcs}
\suca^{1/4}\sqrt\loge< \bfcu\suca^{1/4}\,,
$$
for $E,t$ satisfying \eqref{ninetto}.
Consequently, by \eqref{schiaffino},  the function
$x_1(E,t)$ defined in \eqref{enna2}
is holomorphic\footnote{We are considering the square root as a holomorphic
function in the complex plane excluding the negative real axis.}
for $E,t$ satisfying \eqref{ninetto},
taking $\bfcs$ small enough.
\eproof

\nl
We finally pass to action angle variables
defined in \eqref{alice}.
First we observe that 
by \eqref{LEGOk} and \eqref{ciofecak}, 
for $E$ as in \eqref{ninetto} we get
 \beq{uvaspina}
 |\log \loge|/\sqrt\suca\lessdot
 |\partial_E \act|\lessdot |\log \loge|/\sqrt\suca\,.
 \eeq 
Note that, by
 \eqref{arista} and \eqref{ferrari}, 
 $a_+(\loge)=\act(\tilde E)$.
 By \eqref{uvaspina}
 the image of the ball
 $|E-\tilde E|\leq\bfcs\suca\loge$
 through the function $I(E)$
 contains the ball $|I-a_+(\loge)|\leq \rho$
 (defined in \eqref{blueeyes})
 taking $\hcc$ large enough.
 Analogously
  by \eqref{uvaspina},
  for every 
  $|E-\tilde E|\leq\bfcs\suca\loge$,
 the image of the ball
 $|t|<\bfcs /\sqrt\suca$
 through the function 
 $t\to t/\partial_E \act(E)$
 contains the ball $|\f|\leq 
\s$
 taking $\hcc$ large enough.
Recalling \eqref{playmobilk}, this  complete the proof\footnote{Close to a hyperbolic point
the estimates for the action analyticity radius in \eqref{blueeyes}
is
$\rho	= \loge|\log \loge| \sqrt\suca/C$
since $\partial_E \act\sim |\log \loge|/\sqrt\suca$ (see \eqref{uvaspina}), being $\loge\suca$ the distance in energy
from the critical energy of the hyperbolic point (see \eqref{topo} below).
Far away from the hyperbolic point the derivative is smaller
 (namely, $\partial_E \act\sim 1/\sqrt\suca$) but  the distance in energy is bigger (being $\sim \suca$).}
 of
 \eqref{blueeyes}.

\bigskip

\nl
Let us now prove  \eqref{ofena}. First
by the chain rule we get
\begin{eqnarray}
&&\partial_{\act_1} {\mathtt E}  
=
\frac{1}{\partial_E \act_1 }\,,\qquad
\partial_{\hat \act} {\mathtt E}  
=
-\frac{\partial_{\hat \act} \act_1 }{\partial_E \act_1 }\,,
\qquad
\partial^2_{\act_1\act_1} {\mathtt E}  
=
-\frac{\partial^2_{EE} \act_1 }{(\partial_E \act_1 )^3}\,,
\nonumber
\\
&&\partial_{\act_1 \hat \act}^2 {\mathtt E}  
=
\frac{\partial^2_{EE} \act_1  \partial_{\hat \act}\act_1 }{(\partial_E \act_1 )^3}-
\frac{\partial_{E\hat \act}^2 \act_1  }{(\partial_E \act_1 )^2}
\,,
\label{daitarn3}
\\
&&\partial^2_{\hat \act\hat \act} {\mathtt E}  
=
-\frac{\partial^2_{\hat \act\hat \act} \act_1 }{\partial_E \act_1 }
+\frac{\partial_{\hat \act}^T\act_1\ \partial_{\hat \act}(\partial_{E} \act_1)
+ \partial_{\hat \act}^T(\partial_{E} \act_1)\ \partial_{\hat \act}\act_1 }
{(\partial_E \act_1 )^2}
-\frac{\partial^2_{EE} \act_1 \  \partial_{\hat \act}^T \act_1\  \partial_{\hat \act}\act_1}{(\partial_E \act_1 )^3}
\,,
\nonumber
\end{eqnarray}
where the derivatives of ${\mathtt E} $ and  $\act_1$ are evaluated in  $\big(\act_1 (E,\hat \act),\hat \act\big)$
and $(E,\hat \act),$ respectively.

\nl
 Then, we split $\big(\Bu (\loge)\big)_{\!\rho}$ in two subsets: the region where $\act_1$ is close to
$\acci_+(\loge)$
(near the hyperbolic equilibrium)
and the 
region far away from it.
More precisely, recalling  \eqref{autunno2},
we set
\beqa{sarago}
{\mathcal E}_{\rm cl} &:=&
\Big\{|E-E_+|<\suca/\cc\ |\ \Im(E-E_+)=0 \implies 
\Re(E-E_+)>0\Big\}
\times \hat D_{\ro/4}
\,,
\nonumber
\\
{\mathcal E}_{\rm aw}
&:=&
\mathcal E_{1/4\cc}\cap\{|\Im E|<\suca/(3\cc)^5\}\times \hat D_{\ro/4}\,.
\eeqa
Then,
$$
\big(\Bu (\loge)\big)_{\!\rho}\, \subset\,
{\mathcal E}_{\rm cl}\cap{\mathcal E}_{\rm aw}\,,
$$
taking $\hcc$ large enough.
Regarding the first region we start noting that by 
Cauchy estimates, \eqref{LEGOk} and
\eqref{pappagallok},
for
$
\act_1=\act_1\big(E_\pm(\hat \act)\mp \suca z, \,\hat \act\big)
$, with $|z|<1/\cc$ not belonging to the negative real semiaxe
and $ \hat \act\in \hat D_{\ro/4}$, we have 
\begin{eqnarray}
&&
\frac{|\ln z|}{\sqrt\suca}
\lessdot
|\partial_E\act_1|\lessdot \frac{|\ln z|}{\sqrt\suca}\,,\qquad
\frac{1}{\suca^{3/2}|z|}
\lessdot
|\partial_{EE}\act_1|\lessdot \frac{1}{\suca^{3/2}|z|}\,,
\nonumber
\\
&&
|\partial_{\hat\act}\act_1|\lessdot \lella\,,\qquad
|\partial^2_{E,\hat\act}\act_1|\lessdot \lella|\ln z|\,,\qquad
|\partial^2_{\hat\act,\hat\act}\act_1|
\lessdot\frac{\lella}{\ro}\,,
\label{topo}
\end{eqnarray}
where the first line follows by \eqref{ciofecak}.
Then 
\eqref{ofena} directly follows from
\eqref{daitarn3}.
\\
Consider now the second region, 
namely
 ${\mathcal E}_{\rm aw}$.
 By
 \eqref{rosettaTH} with
 $\l=1/\cc$
  we have
 $|\partial_E\act_1|\leq 5\cc^3/\sqrt\suca$
  on
 $\mathcal E_{1/2\cc}\times \hat D_{\ro/2}$.
 Then by Cauchy estimates
 we get
 $|\partial^2_{EE}\act_1|\leq 40\cc^4/\suca^{3/2}$
 on ${\mathcal E}_{\rm aw}$.
 Then by \eqref{vana} we get
 $|\partial_E\act_1|\geq 1/2\cc\sqrt\suca$
 on ${\mathcal E}_{\rm aw}$.
 Using this lower bound, \eqref{daitarn3},
 \eqref{rosettaTH}  and Cauchy estimates,
 \eqref{ofena} follows also in 
 ${\mathcal E}_{\rm aw}$. 
 
 \medskip
 
\nl
We finally prove
\eqref{inthecourtk}.
Since the maps $\cFiq$ preserve the $(n+1)$-dimensional measure
$dId\f_1$, recalling \eqref{arista}--\eqref{playmobilk}
and using 
\eqref{LEGOk}--\eqref{pappagallok},
\eqref{cimabue}, \eqref{alce} one obtains, by Fubini's theorem,
\begin{eqnarray*}
&&\meas\Big(
\big(D^\flat\times\T\big)\ \setminus\ 
\bigcup_{0\leq i\leq 2N}\check\cM^i (\loge)\Big)
\leq 
2\pi\sum_{0\leq i\leq 2N} \meas(\hat D) \meas\left( \Bu^i (0)
\setminus\Bu^i (\loge)\right)
\\
&&\leq 
2\pi\sum_{0\leq i\leq 2N} \meas(\hat D)
\sup_{\hat\act\in\hat D}
\left(
\acci_+(\hat \act;0)
-
\acci_+(\hat \act;\loge)
+
\acci_-(\hat \act;\loge)
-
\acci_-(\hat \act;0)
\right)
\\
&&\stackrel{\eqref{rosettaTH}}\lessdot 
 \meas(\hat D)\sqrt\suca
\int_0^\loge
|\log z| \, dz
\leq
\meas(\hat D) \sqrt\suca \loge |\log\loge|\,. \qedeq
\end{eqnarray*}


\section{Convexity energy estimates}\label{convexity}

In this section we investigate the convexity of the energy functions $\act_1\to {\mathtt E}^{i}(\act_1,\hat{p})$ defined as the inverse functions of the action functions\footnote{Observe that  the action function $E\to\act_1^{i}(E,\hat{p})$ is strictly increasing and hence invertible.} $E\to\act_1^{i}(E,\hat{p})$. \\
We also denote 
$\bar\act_1^{i}:=\act_1^{i}|_{\lalla=0}$ the `unperturbed action function'  and its
inverse $\bar{\mathtt E}^{i}:=
{\mathtt E}^{i}|_{\lalla=0}$ the `unperturbed energy function'.

\rem\label{maschio}  Observe that 
$\bar\act_1^{(0)}(E)=\bar\act_1^{(2N)}(E)$
and $\bar{\mathtt E}^{(0)}(\act_1)=\bar{\mathtt E}^{(2N)}(\act_1)$.
\erem

\nl
In general, the energy functions have inflection points\footnote{Compare \cite{BCsecondary}.}, 
however there are some cases in which the convexity is definite, namely, in the outer regions ($i=0,2N$) and in the case the 
reference potential $\GO$ is `close' to a cosine  (in which case $N=2$) in the sense of the following

\dfn{pigro}{\bf (Cosine--like functions)}
 Let  $0<\ttg< 1/4$.
 We say that a real analytic function $G:\T_1\to\C$  is
$\ttg$--cosine--like
if, for  some $\eta>0$ and 
 $\sa_0\in\R$, one~has
  \begin{equation}
 \sup_{\sa\in\torus_1}{\modulo}G(\sa)-\eta\cos (\sa+\sa_0) {\modulo}
 \leq \eta\ttg\,.
 \label{A3bis}
\end{equation}
 \edfn

\begin{proposition}\label{kalevala} {\rm (i)}
If $i=0,2N$, then, for every $E> \bar E_{i}$, one has: 
$\partial^2_{\act_1}\bar{\mathtt E}^{i}(\bar \act_1^{i}(E))
\geq 2$.

\nl
{\rm (ii)}
If $\GO$ is cosine--like with $\ttg\leq 2^{-40}$,
then
\begin{equation}\label{thor}
\partial^2_{\act_1} \bar{\mathtt E}^{1}           
(\bar \act_1^{1}(E))
\, \leq \,
-\frac{1}{27}\,,
\qquad
\forall E\in(\bar E_1,\bar E_2)
\,.
\end{equation}
\end{proposition}

 \proof   
 (i) 
 Let us consider now the zone above separatrices.
 First observe that the cases $i=0$ and 
$i=2N$  are identical by Remark \ref{maschio}. 
Let us then consider the case $i=2N$. By definition,
\begin{equation*}
 \bar \act_{1}^{2N}(E)
=
\frac{1}{2\pi}\int_0^{2\pi}\sqrt{E-\GO(x)}dx\,,
 \end{equation*}
 thus, by   Jensen's inequality
 $$
 (2\partial_E\bar \act_1^{i}(E))^3=
 \Big(\frac{1}{2\pi}\int_0^{2\pi}\frac1{\sqrt{E-\GO(x)}}dx\Big)^3
 \leq 
 \frac{1}{2\pi}\int_0^{2\pi}\frac1{(E-\GO(x))^{3/2}}dx
 =-4 \partial^2_{E}\bar \act_1^{i}(E)\,,
 $$
 and the claim follows by
 \begin{equation}\label{doraemon}
\partial^2_{\act_1 \act_1}\bar{\mathtt E}^{(i)}\big(\bar\act^{(i)}_1(E)\big)
=-\frac{\partial_{EE}^2\bar\act^{(i)}_1(E)}{\big(\partial_{E}\bar\act^{(i)}_1(E)\big)^3}\,.
\end{equation} 
(ii)  
First we note that, up to a phase translation,  we can take $\sa_0=0$
in \eqref{A3bis}.
Then set
\begin{equation}\label{barbacane}
M:=\max_{\R} \GO\,,\quad
m:=\min_{\R} \GO\,,\quad
L(y):=\frac{2y-M-m}{M-m}\,,\quad
V:=L\circ \GO\,.
\end{equation}
Note that $\max_{\R} V=1,$ $\min_{\R}V=-1$.\\
The idea is to study the action variable of the Hamiltonian 
$p_1^2+V(q_1)$ which is strictly related the one of $p_1^2+\GO(q_1)$,
see \eqref{caliburnus} below.
Denoting $|\cdot|_r:=\sup_{\torus_r}|\cdot|$, we have the following

\begin{lemma}\label{fossato}
If $\GO$ satisfies \eqref{A3bis} then $V$ in \equ{barbacane} satisfies 
$|V-\cos z|_1\leq 4\ttg$.
\end{lemma}
\proof
By  \eqref{A3bis} and \eqref{barbacane} we have
\begin{equation}\label{postierla}
-\ttg\ch\,\leq \, M-\ch\,,\ m+\ch\,\leq\, \ttg\ch\,,\qquad
\left|\frac{2\ch}{M-m}-1\right|\leq \frac{\ttg}{1-\ttg}\,.
\end{equation}
By  \eqref{A3bis} and \eqref{barbacane} we get
$$
\left|V(z)-\frac{2\ch}{M-m}\cos z\right|_1\leq
\frac{2}{M-m}\left(
\ttg\ch+\frac{M+m}{2}
\right)
\leq
\frac{4\ttg\ch}{M-m}\leq\frac{2\ttg}{1-\ttg}\,.
$$
Then
$$
\sup_{\torus_1}\left|V(z)-\cos z\right|\leq
\frac{2\ttg}{1-\ttg}
+
\left|\frac{2\ch-M+m}{M-m}\cos z\right|_1
\leq
\frac{2+\cosh 1}{1-\ttg}\leq 4\ttg
\,.\qedeq
$$
Next we need a representation lemma whose proof is given in Appendix~\ref{stockausen}:
\begin{lemma}\label{parcheggio} Let $0<\tgo\le 2^{-10}$ and let
${w}$ be a real analytic $2\pi$-periodic function  satisfying
 \beqno
 \max_\R {w}=1\,,\quad \min_\R {w}=-1\,,\quad
 |{w}(z)-\cos z|_1\leq \tgo \,.
 \eeqno
Then,  there exists a unique 
 real analytic $2\pi$--periodic function $b$ such that  such that 
 $${w}(z)=\cos(z+b(z))\,,\qquad |b|_{1/4}\leq 9\sqrt\tgo\,.
 $$
\end{lemma}

\noindent
Lemma \ref{parcheggio} can be applied to  the potential $V$ in \equ{barbacane}, so that, in particular,  $V$ has only two critical points
(a maximum and a minimum) on a period.
\\
For $i=0,1,2$ let us denote by $\bar\act^{(i)}_1(E)$, respectively, $\tilde\act^{(i)}_1(E),$ the action variable 
of the Hamiltonian $p_1^2+\GO(q_1)$, respect.  $p_1^2+V(q_1)$
in the three zones below ($i=0$), inside ($i=1$) and above ($i=2$)
separatrices.

\nl
Now, the relation between the action  $\bar\act^{(i)}_1$ of $p_1^2+\GO(q_1)$ and the action $\tilde\act^{(i)}_1$ of $p_1^2+V(q_1)$  is given by the following formula:
\begin{equation}\label{caliburnus}
\bar\act^{(i)}_1(E)=\sqrt{\frac{M-m}{2}}\tilde\act^{(i)}_1\big(L(E)\big)\,,\qquad i=0,1,2\,.
\end{equation}
Indeed, considering the case $i=2$ (the other ones being analogous), and recalling \eqref{barbacane}, one finds
\begin{eqnarray*}
\bar\act^{(i)}_1(E)
&=&\frac{1}{2\pi}\int_0^{2\pi}\sqrt{E-\GO(x)}dx
=
\frac{1}{2\pi}\int_0^{2\pi}\sqrt{L^{-1}(L(E))-L^{-1}(V)(x)}dx
\\
&=&\sqrt{\frac{M-m}{2}}
\frac{1}{2\pi}\int_0^{2\pi}\sqrt{L(E)-V(x)}dx
=\sqrt{\frac{M-m}{2}}\tilde\act^{(i)}_1\big(L(E)\big)\,,
\end{eqnarray*}
which proves \equ{caliburnus}.

\nl
Going back to the proof of \equ{thor},
denote by 
$\bar{\mathtt E}^{(i)}(\act_1)$, respectively $\tilde{\mathtt E}^{(i)}(\act_1)$,
the inverse function of
$\bar\act^{(i)}_1(E)$, respectively $\tilde\act^{(i)}_1(E)$.
By\footnote{And the  analogous formula for 
$\partial^2_{\act_1 \act_1}\tilde{\mathtt E}^{(i)}\big(\tilde\act^{(i)}_1(E)\big)$.} \equ{doraemon}, \eqref{caliburnus} and \eqref{barbacane} we get
 \begin{equation}\label{doraemon2}
\partial^2_{\act_1 \act_1}\bar{\mathtt E}^{(i)}\big(\bar\act^{(i)}_1(E)\big)
=
\partial^2_{\act_1 \act_1}\tilde{\mathtt E}^{(i)}\big(\tilde\act^{(i)}_1(L(E))\big)
\end{equation}
Let us consider first the zone inside separatrices and, to simplify notation, 
denote   $\tilde\act_1^{(1)}(E)$
by $A(E)$. Then,
$$
A(E)=A_-(E)+A_+(E):=
\frac{1}{\pi}\int_{x_m-2\pi}^{V_-^{-1}(E)}
\sqrt{E-V(x)}dx
+
\frac{1}{\pi}\int_{V_+^{-1}(E)}^{x_m}
\sqrt{E-V(x)}dx\,,
$$
where $V_-^{-1}(E)$ and $V_+^{-1}(E)$ are the inverse of 
$V(x)=\cos(\psi(x))=E$, with $\psi(x):=x+b(x)$,
in the intervals $[x_m-2\pi,x_M]$ and $[x_M,x_m],$
respectively; namely
$V_-^{-1}(E)=\psi^{-1}(-\arccos(E))$ and
$V_+^{-1}(E)=\psi^{-1}(\arccos(E))$.
Recall that $\psi(x_M)=0$ and $\psi(x_m)=\pi$.
Since $|b|_{1/4}\leq 18\sqrt\ttg$ by Cauchy estimates
we have that $\psi$ is invertible with inverse
$\psi^{-1}(y)=y+u(y)$ for a suitable\footnote{$u$
is the solution of the fixed point equation $u(y)=-b(y+u(y))$
in the space of $2\pi$-periodic real 
analytic function with holomorphic extension on the strip
$\{|\Im y|<1/5\}$ and $|u|_{1/5}\leq 18\sqrt\ttg$.} 
$2\pi$-periodic real 
analytic
$u$ satisfying 
\begin{equation}\label{briscola}
|u|_{1/5}\leq 18\sqrt\ttg.
\end{equation}
We get
$$
A_+'(E)=
\frac{1}{2\pi}\int_{V_+^{-1}(E)}^{x_m}
\frac{dx}{\sqrt{E-V(x)}}=
\frac{1}{2\pi}\int_0^1
\frac{1+u'\big(\arccos(g(E,t))\big)}{\sqrt{t-t^2}
\sqrt{1+t-E(1-t)}}dt\,,
$$
making the substitution $x=x(t):=V_+^{-1}(g(E,t))$ with
$g(E,t):=E-(1+E)t$.
Analogously
$$
A_-'(E)=
\frac{1}{2\pi}\int_0^1
\frac{1+u'\big(-\arccos(g(E,t))\big)}{\sqrt{t-t^2}
\sqrt{1+t-E(1-t)}}dt\,.
$$
Then taking the even part $v$ of $u'$, namely
$v(y):=\frac12(u'(y)+u'(-y))$
we have
$$
A'(E)=
\frac{1}{\pi}\int_0^1
\frac{1+v\big(\arccos(g(E,t))\big)}{\sqrt{t-t^2}
\sqrt{1+t-E(1-t)}}dt\,.
$$
Note that by Cauchy estimates 
$|v|_{1/6}\leq 540\sqrt{\ttg}$.
Deriving we get
$$
A''(E)=
\frac{1}{2\pi}\int_0^1\sqrt{\frac{1-t}{t}}\,
\frac{1+v_0(E,t)}{(1+t-E(1-t))^{3/2}}\,
dt\,,
$$
with
\begin{eqnarray*}
v_0(E,t)&:=&
v\big(\arccos(g(E,t))-2\tilde v(E,t)\,,
\\
\tilde v(E,t)
&:=&
\frac{v'\big(\arccos(g(E,t))\big)\, \sqrt{1+t-E(1-t)}}{\sqrt{1-t}\,
\sqrt{1+E}}\,.
\end{eqnarray*}
Since $v$ is $2\pi$-periodic and even, we have $v'(\pi)=0$.
Then, by Cauchy estimates we get
$$|v'(\xi)|\leq 39880\cdot\sqrt{\ttg}\  |\xi-\pi|\,,\qquad \forall\ \xi\in\R\,.
$$
Note that 
$$0\leq \pi-\arccos(-1+\xi)\leq \frac{\pi}{\sqrt 2}\sqrt{\xi}\,,\qquad \forall\ 0\leq\xi\leq 2\,.
$$
Therefore, since $g(E,t)+1=(1-t)(1+E)$, 
for  $0<t<1$ and $-1<E<1$, one has
$$
|v'\big(\arccos(g(E,t))\big)|\leq 39880\sqrt\ttg
|\pi-\arccos(g(E,t))|
\leq 19440\pi\sqrt {2 \ttg} \sqrt{1-t}\,
\sqrt{1+E}\,,
$$
which implies 
$$
|\tilde v(E,t)|\leq 244292\sqrt\ttg\qquad
{\rm and }\qquad
|v_0(E,t)|\leq 
244292\sqrt\ttg\leq 
2^{18}\sqrt\ttg\,.
$$
Taking
$\ttg\leq 2^{-38}$
we have $|v_0(E,t)|\leq 1/2$ and therefore
for every $-1<E<1$
\begin{equation}\label{giada}
\frac12 A_0'(E)\leq A'(E)\leq \frac32 A_0'(E)\,,\qquad
\frac12 A_0''(E)\leq A''(E)\leq \frac32 A_0''(E)\,,
\end{equation}
where $A_0(E)$ is the action variable with exactly
cosine potential (namely when $\ttg=0$), namely
\begin{eqnarray*}
A_0'(E)
&=&
\frac{1}{\pi}\int_0^1
\frac{1}{\sqrt{t-t^2}
\sqrt{1+t-E(1-t)}}dt\,,
\\
A_0''(E)
&=&
\frac{1}{2\pi}\int_0^1\sqrt{\frac{1-t}{t}}\,
\frac{1}{(1+t-E(1-t))^{3/2}}\,
dt\,.
\end{eqnarray*}
Then, since $\tilde{\mathtt E}^{(1)}(I_1)$ is the inverse of 
$\tilde\act_1^{(1)}(E)=A(E)$,
 for every $-1<E<1$
\begin{equation}\label{nobita}
-\partial^2_{\act_1} \tilde{\mathtt E}^{(1)}
(\tilde \act_1^{(1)}(E))=
\frac{A''(E)}{(A'(E))^3}\geq \frac{4}{27}\frac{A_0''(E)}{(A_0'(E))^3}\geq \frac{1}{27}\,,
\end{equation}
since, as it not difficult to check,  the function $\frac{A_0''(E)}{(A_0'(E))^3}$ is increasing and has limit $1/4$
for  $E\to-1^+$.
\qed


\appendix

\section{Proofs of Proposition~\ref{bruegel}}\label{oceania}

First, recalling \eqref{beatoangelico} and 
\eqref{caspiterina},
we define the symplectic transformation
\beq{ilvedovo}
 \Phi_*: (\ttp,\ttq)\in D_{7\ro/8,7\so/8}
 \ \longrightarrow\
 \big(\ttp,\ttq_1+\sa_{2j}(\hat \ttp),\hat \ttq+\ttp_1\partial_{\hat \ttp}
 \sa_{2j}(\hat \ttp)\big)\in
 D_{\ro,\so}\,,
\eeq
 transforming the Hamiltonian $\Hpend$ in \equ{pasqua}
 into
\beq{toledo}
 \Hpend_*:=\Hpend\circ\Phi_*(\ttp,\ttq)=:
 \big(1+ \cin_*(\ttp,\ttq_1)\big) \ttp_1^2  
  +\Gm_*(\hat \ttp, \ttq_1)\,.
\eeq
By Taylor expansion at
$(\ttp,\ttq_1)=(0,\hat \ttp,0)$,
recalling \equ{toledo}, \equ{cimabue2}, \equ{pinturicchio} and \equ{alce},
we get
\beqa{toledo2}
 &&\Hpend_*=E_{2j}(\hat \ttp)+
 \big(1+ \cin_*(0,\hat \ttp,0)\big) \ttp_1^2  
  -\l^2(\hat \ttp)\ttq_1^2+R_*(\ttp,\ttq_1)
  \,,\qquad \mbox{with}
  \nonumber
  \\
  &&R_*(\ttp,\ttq_1):=\big(\cin_*(\ttp,\ttq_1)-\cin_*(0,\hat \ttp,0)\big) \ttp_1^2
  +\Gm_*(\hat \ttp, \ttq_1)-
  \frac12\partial^2_{\ttq_1}\Gm_*(\hat \ttp, 0)\ttq_1^2
  \,.
\eeqa
Then, the following trivial lemma holds.
\begin{lemma}
  There exist a constant $0<\bfco<1/8$, depending only on 
 $\upkappa,n$, such that, defining 
 the symplectic transformation
\begin{eqnarray}\label{lavedova0}
&& \Phi_0:  \{|Y_1|<\bfco \suca^{1/4}\}\times\hat D_{3\ro/4}
 \times \{|X_1|<\bfco \suca^{1/4}\}\times \T^{n-1}_{3\so/4}
 \ \longrightarrow\
 D_{7\ro/8,7\so/8}\,,
\\
&&
 \ttp_1=\d(\hat Y)Y_1,
 \qquad
 \hat \ttp=\hat Y,
\qquad
 \ttq_1=\frac1{\d(\hat Y)} X_1\,,
 \qquad
\hat \ttq=\hat X-\frac{\partial_{\hat Y}\d(\hat Y)}{\d(\hat Y)}Y_1 X_1\,,
\nonumber
\end{eqnarray}
 we have that $ \Hpend_0:=	
 \Hpend_*\circ \Phi_0$ has the form
 \beq{avila}
 \Hpend_0=E_{2j}(\hat Y)
 +g(\hat Y)(Y_1^2-X_1^2)+\suca
 R_0(\suca^{-1/4}Y_1,\hat Y,\suca^{-1/4}X_1)
\,,
	\eeq
	where $R_0(\tilde Y_1,\hat Y,\tilde X_1)$ is holomorphic on 
	$$
	\{|\tilde Y_1|<\bfco \}\times\hat D_{3\ro/4}
 \times \{|\tilde X_1|<\bfco \}\,,
	$$
	with $|R_0|\lessdot 1$
	and, finally,
	it is at least cubic in $\tilde Y_1,\tilde X_1$.
\end{lemma}
\proof The fact that $\Phi_0$ is well defined on its domain
follows by the explicit expression in 
\equ{lavedova0}, by \equ{deltoide}, \equ{alce} and \equ{cimabue} (in particular $\suca\leq \ro^2/2^{16}$).
Eq. \equ{avila} follows by \equ{toledo2} setting
\beq{toledo3}
R_0(\tilde Y_1,\hat Y,\tilde X_1):=\suca^{-1}
 R_*\left(\d(\hat Y)\suca^{1/4}\tilde Y_1,\,\hat Y,\,
\frac{\suca^{1/4}}{\d(\hat Y)}\tilde X_1\right)\,.
\eeq
Finally,  the estimate $|R_0|\lessdot 1$  follows from \equ{cimabue2}
and \equ{deltoide}.
\eproof
Next, we shall use the following well known result, whose proof can be found, e.g., in\footnote{See, in particular, Lemma 0 and Appendix A.3 in \cite{CG}.}
\cite{CG} or in \cite{G}. 

\begin{lemma}\label{redhair}
Given a Hamiltonian $\Hpend_0$ as in
\equ{avila}.
For suitable constants $0<\bfcu<\bfco/8n\bfcd$,  depending only on 
 $\upkappa,n$, 
 there exist
 a (close to the identity) symplectic transformation
 \begin{eqnarray}\label{lavedovabis}
 \Phi_1:  &&\{|y_1|<\bfcu \suca^{1/4}\}\times\hat D_{\ro/2}
 \times \{|x_1|<\bfcu \suca^{1/4}\}\times \T^{n-1}_{\so/2}
 \ \longrightarrow
 \\
&&\{|Y_1|<\bfco \suca^{1/4}\}\times\hat D_{3\ro/4}
 \times \{|X_1|<\bfco \suca^{1/4}\}\times \T^{n-1}_{3\so/4}\,,
 \nonumber
 \end{eqnarray}
and  a function $R_{\rm hp}(z,\hat y)$ 
satisfying \eqref{molecolare} and \eqref{biretta},
 such that
 $ \Hpend_{\rm hp}(y,x_1):=\Hpend_0\circ \Phi_1(y,x)$
 satisfies \eqref{abbacchio}.
Moreover $\Phi_1$ has the form
 \beqa{pesantebis}
 &&Y_1= y_1+\suca^{1/4} a_1(\suca^{-1/4}y_1,\hat y,\suca^{-1/4}x_1)\,,\quad 
 \hat Y=\hat y,\\
 &&X_1=x_1+
 \suca^{1/4} a_2(\suca^{-1/4}y_1,\hat y,\suca^{-1/4}x_1),\quad 
 \hat X=\hat x+
 \sqrt\suca \ro^{-1} 
 a_3(\suca^{-1/4}y_1,\hat y,\suca^{-1/4}x_1)\,,
\nonumber
 \eeqa
for suitable functions $a_i(\tilde y_1,\hat y,\tilde x_1)$, $i=1,2,3$,
which are 
 holomorphic and bounded by $\bfcd$ on 
	$$
	\{|\tilde y_1|<\bfco/2 \}\times\hat D_{\ro/2}
 \times \{|\tilde x_1|<\bfco/2 \}\,,
	$$
	moreover, $a_1,a_2$, respectively
	$a_3$, are
	at least quadratic, respectively cubic, in $\tilde y_1,\tilde x_1$.
\end{lemma}

\rem\label{domiziano}
$ \Hpend_{\rm hp}$ is simply the hyperbolic Birkhoff normal form of $\Hpend_0$.
Any canonical transformation of the form
$ y_1=\a\tilde y_1+\b\tilde x_1,$ $ x_1=\b\tilde y_1+\a\tilde x_1,$ with $\a^2-\b^2=1$
and $\hat y=\hat{\tilde y}$ leaves $ \Hpend_{\rm hp}$ invariant since
$y_1^2-x_1^2=\tilde y_1^2-\tilde x_1^2.$
Namely the integrating transformation $\Phi_1$ is not unique.
However, as well known,  the form of the integrated Hamiltonian 
$\Hpend_{\rm hp}$ in \equ{abbacchio} is unique,
in the sense that $E_{2j}$, $g$ and $R$ are unique.
\\
Note also that
the map $\Phi_1$ is close to the identity,
for $\bfcu$ small,
since  its Jacobian  
is the identity plus a matrix whose entries are
 (by Cauchy estimates) uniformly bounded
 on its domain in \equ{lavedova} by 
 $2\bfcd\bfcu/\bfco\leq 1/4n$.
\erem
Leu us go back to the proof of  Proposition~\ref{bruegel}  and let us   prove  
 \eqref{biretta}.
 \\
 Evaluating \equ{abbacchio} for $\lalla=0$ we get
\beqa{abbacchio2}\nonumber
\bHpend_1(y,x_1)&:=&\bHpend_1(y,x_1)|_{\lalla=0}=\bHpend_0\circ \bar\Phi_1(y,x)\\
&=& \bar E_{2j}+\bar g (y_1^2-x_1^2)
+\suca \bar R_{\rm hp}\left(\frac{y_1^2-x_1^2}{\sqrt\suca} \right)=O(\suca)
\eeqa
on the domain defined in \equ{lavedova}.
Let us denote
$\bHpend_0:=\Hpend_0|_{\lalla=0}$.
Since
by \eqref{toledo2}, \eqref{toledo3}, \eqref{cimabue2}
one has $\Hpend_0-\bHpend_0=O(\suca\lalla)$,
$$
\Hpend_0\circ \bar\Phi_1=\bHpend_0\circ \bar\Phi_1
+(\Hpend_0-\bHpend_0)\circ \bar\Phi_1=
\bHpend_1 +R_1\,,\quad \mbox{with}\ \ \ 
R_1=O(\suca\lalla)\,,
$$
namely the system is integrated up to a small term of order
$\suca\lalla$.
Note also that, since $\bar\Phi_1$ has the form in \equ{pesantebis},
it leaves invariant the terms of order $\leq 2$ in $(y_1,x_1)$, namely
\beq{maestra}
\Hpend_0\circ \bar\Phi_1=E_{2j}+g(y_1^2-x_1^2)+\suca\bar R+Q\,, 
\quad \mbox{with} \quad Q=O(\suca\lalla)\,.
\eeq
 Now we want to construct a symplectic transformation
$\Phi_\lalla$ integrating $\Hpend_0\circ \bar\Phi_1$.
Since $\bHpend_1$ is already in normal form,
we claim that
 the integrating transformation $\Phi_\lalla$ is $O(\suca^{1/4}\lalla)$--close to the identity
and 
\beq{nottefonda}
\Hpend_0\circ \bar\Phi_1\circ\Phi_\lalla=(\bHpend_1 +R_1)\circ\Phi_\lalla
=:\Hpend_{\rm hp}'=\bHpend_1 +O(\suca\lalla),
\eeq
where $\Hpend_{\rm hp}'$ is in normal form, namely as the form in 
\equ{abbacchio}.
By the unicity of the Birkhoff Normal Form
we deduce that $\Hpend_{\rm hp}=\Hpend_{\rm hp}'=\bHpend_1 +O(\suca\lalla)$.
By \eqref{abbacchio}, \eqref{abbacchio2}, \eqref{deltoide}
and \eqref{october} we get \equ{biretta}.

\nl
It remains to prove
\equ{nottefonda}.
The crucial 
 point here is that the generating 
 function\footnote{According to the Lie's series method}
 $\chi$
 of the integrating transformation $\Phi_\lalla$  is 
$O(\sqrt\suca\lalla)$ and its gradient is, by Cauchy estimates,
$O(\suca^{1/4}\lalla)$ in a domain
$\{|y_1|,|x_1|\lessdot\suca^{1/4}\}$.
The fact that $\chi=O(\sqrt\suca\lalla)$ can be easily seen passing, as usual in Birkhoff Normal Form, to the coordinate $\xi=(y_1-x_1)/\sqrt 2$,
$\eta=(y_1+x_1)/\sqrt 2$. In these coordinates, recalling 
\equ{maestra}, we get
$$
\Hpend_0\circ \bar\Phi_1=E_{2j}+2g\xi\eta+ \suca \bar R'(\xi,\eta)+Q'(\xi,\eta)$$
with 
$ \bar R'=\bar R_{\rm hp}(2\xi\eta/\sqrt\suca)=O(1)$ and $Q'=O(\suca\lalla)
$.
Note that
 the Taylor expansion of
$\bar R'$ containing only monomial of the form $\bar R'_{hh} \xi^h \eta^h$.
At the first step, we have to cancel all the monomials of $Q'$
of the form $Q'_{hk} \xi^h \eta^k$ with $h+k=3$.
The generating function$\chi^{(3)}$
of the first step is exactly
$$
\chi^{(3)}=\sum_{h+k=3}\frac{Q'_{hk}}{2g(h-k)}\xi^h \eta^k\stackrel{\equ{deltoide}}=O(\sqrt\suca\lalla)\,.
$$
After this first step the Hamiltonian becomes
$
E_{2j}+2g\xi\eta+ \suca \bar R'(\xi,\eta)+Q''(\xi,\eta)$ with  $Q''=O(\suca\lalla)
$.
At the second step, we have to cancel all the monomials of $Q''$
of the form $Q''_{hk} \xi^h \eta^k$ with $h+k=4$, $h\neq k$. 
We proceed as in the first step with analogous estimates.
Analogously for the other infinite steps, obtaining
$$
E_{2j}+2g\xi\eta+ \suca \bar R'(\xi,\eta)+\bar Q(\xi,\eta)$$ with  $\bar Q=O(\suca\lalla)$ and $\bar Q_{hk}=0$ for $h\neq k$,
proving \equ{nottefonda} (recall \eqref{deltoide}
and \eqref{october}).
\eproof
We can conclude the proof of Proposition
\ref{bruegel}:\\
 The composition of the symplectic transformations
defined in \equ{ilvedovo}, \equ{lavedova0}, \equ{lavedova}
integrates $\Hpend$, namely
\eqref{abbacchio} holds\footnote{As well as
\eqref{molecolare} and \eqref{biretta} by Lemma 
\ref{redhair}.}  with
$\Phi_{\rm hp}:=\Phi_*\circ\Phi_0\circ\Phi_1$
satisfying \eqref{lavedova}, \eqref{pesante}
and \eqref{chedolore}. 
\\
The inclusion \eqref{brodino}
follows by \eqref{pesante} and \eqref{deltoide}.
\eproof

\section{Proofs of two simple lemmata}\label{australia}

\subsection{Proof of Lemma~\ref{enricone}} \label{cento}

We know that $\partial_\sa \GO(\bar\sa_i)=0 $
and we want to solve 
the  equation
$
\partial_\sa \Gm(\hat{p},\sa_i(\hat{p}))=0$.
Equivalently, for $\lalla\leq 2^{-8}\upkappa^{-6}$, we want to find a real analytic
$y=y(\hat{p}),$ $\hat{p} \in \hat D_\ro$,
with 
\beq{falanghina}\ts
\sup_{\hat D_\ro}|y|\leq \rho
:= \frac{2\suca\lalla}{\morse \so}
\stackrel{\equ{alce}}\leq \frac{\so}{2}\,,
\eeq
solving the equation
\beq{sterco}
\partial_\sa \Gm(\hat{p} ,\bar\sa_i+y(\hat p))=0\,,
\eeq
so that $\sa_i(\hat p)=\bar\sa_i+y(\hat p)$.
We have
$$
\partial_\sa \Gm(\hat{p} ,\bar\sa_i+y)=
\partial_\sa \Gm(\hat{p} ,\bar\sa_i)+
g(\hat{p},y) y\,,\quad
\mbox{where}
\quad
g(\hat{p},y):=\int_0^1
\partial_{\sa}^2 \Gm(\hat{p} ,\bar\sa_i+ty)dt\,.
$$
Then \equ{sterco} can be written as the fixed point equation
$$
y=\Psi(y)\,,\quad
\mbox{where}
\quad
\Psi(y):=-\frac{\partial_\sa \Gm(\hat{p} ,\bar\sa_i)}{
g(\hat{p},y)}
$$
to be solved in the closed set of the
real analytic functions $y=y(\hat{p})$ on $\hat D_\ro$ satisfying the bound \equ{falanghina}.
Note that, since $\partial_\sa \GO(\bar\sa_i)=0 $, by \eqref{ladispoli} we have 
$|\partial_{\sa}^2 \GO(\bar\sa_i)|\geq \morse$.
Moreover by \equ{cimabue} and Cauchy estimates we get for $|y|\leq \rho$ and 
$\hat{p} \in \hat D_\ro$
\beq{cetriolini}
|g-\partial_{\sa}^2 \GO(\bar\sa_i)|\leq 
\frac{4\suca\lalla}{\so^2}\,, \quad\mbox{which implies}\quad
|g|\geq \morse-\frac{4\suca\lalla}{\so^2}
\stackrel{\equ{alce}}\geq \frac{\morse}{2}\,.
\eeq
Again by $\partial_\sa \GO(\bar\sa_i)=0 $, \equ{cimabue} and Cauchy estimates
we obtain uniformly on $\hat D_\ro$ that
\beq{fragoline}
|\partial_\sa \Gm(\hat{p} ,\bar\sa_i)|\leq
\suca\lalla/\so\,.
\eeq
Then by \equ{cetriolini}
we obtain for $|y|\leq \rho$ and 
$\hat{p} \in \hat D_\ro$
\beq{chiodo}\ts
|\Psi|\leq\frac{2\suca\lalla}{\morse \so}= \rho\,,
\eeq
by \equ{alce} and \equ{falanghina}.
Moreover
$$
\partial_y\Psi(y):=\frac{\partial_\sa \Gm(\hat{p} ,\bar\sa_i)}{
(g(\hat{p},y))^2}\partial_y g(\hat{p},y)\,.
$$
Then for $|y|\leq \rho$ and 
$\hat{p} \in \hat D_\ro$ we get
\beq{chiodo2}\ts
|\partial_y\Psi|< 2^6\frac{\suca^2\lalla}{\morse^2\so^4}
\leq 2^6\upkappa^6\lalla\leq 1\,,
\eeq
by \equ{fragoline}, \equ{cetriolini}, \equ{alce} and since
$
|\partial_y g(\hat{p},y)|<{16\suca}/{\so^3}
$
by \equ{cimabue} and Cauchy estimates.
In conclusion, by \equ{chiodo} and \equ{chiodo2}
we have that $\Psi$ is a contraction and the Fixed Point Theorem applies proving the first estimate
in \equ{october}.
\\
Let us now show the second estimate
in \equ{october}.\\
By \equ{cimabue}, the first estimate
in \equ{october}, \equ{alce} and Cauchy estimates we get
\beqano
|E_i(\hat{p})-\bar E_i|
&\leq& 
|\Gm(\hat p, \sa_j(\hat p))-\GO( \sa_j(\hat p))|
+|\GO( \sa_j(\hat p))-\GO( \bar\sa_j)|
\\
&\leq& 
\ts
\suca\lalla+
\frac{2\suca^2\lalla}{\morse \so^2}
 \leq 3\upkappa^3 \suca\lalla\,,
\eeqano
proving the second estimate in \equ{october}.

\nl
Let us prove the final claim.  By \equ{dicembre} (applied to $\GO$) and by Cauchy estimates it follows that the minimal distance between two critical points of $\GO$ can be estimated from below by $2\morse\so^2/\suca$. Thus, by the first estimates in \equ{october}, it follows that the relative order of the critical points of $\GO$ is preserved, provided
$8 \suca^3\lalla^2<\morse^3\so^4$, which, using \equ{alce} is implied by $2^3 \upkappa^7 \lalla^2<1$, which, in turn, is implied by the hypothesis $\lalla\leq 1/(2\upkappa)^6$. \\
As for critical energies, since $\GO$ is $\morse$--Morse, they are at least $\morse$ apart; hence, from the second estimate in \equ{october} the claim follows provided  $3\upkappa^3 \suca\lalla<\morse$, which by \equ{alce}, is implied by 
$\lalla<1/(3\upkappa^4)$, which, again, is implied by the hypothesis.
\eproof

\subsection{Proof of Lemma~\ref{parcheggio}}\label{stockausen}
First denote $R(z):={w}(z)-\cos z,$ so that $|R|_1\leq \tgo.$
We note that, on the real line, ${w}$ has exactly two critical points: a maximum
 $x_M$ (with ${w}(x_M)=1$) and a minimum $x_m$
 (with ${w}(x_m)=-1$) in the interval
$[-\pi/2,3\pi/2).$ Indeed, since by Cauchy estimates
$\sup_R|{w}'|\leq \tgo,$
 the equation ${w}'(x)=-\sin x+R'(x)=0$ in the interval
$[-\pi/2,3\pi/2)$ has only two solo $x_M,x_m$ with
$|x_M|,|x_m-\pi|\leq 1.0001\tgo\leq 0.001.$ Obviously
$x_M+b(x_M)=0$ and $x_m+b(x_m)=\pi.$
\\
On the real line the function $b$ is given by the $2\pi$-periodic
continuous\footnote{Since $b(x_m-2\pi)=b(x_m)=\pi-x_m$} 
function defined in the interval $[x_m-2\pi,x_m]$ by the expression
$$
b(x):=\sign(x-x_M)\arccos({w}(x))-x
$$
Let us consider first the complex domain
$\Omega_0:=\{0.4< \Re z<\pi-0.4\,,\ |\Im z|<1/4\}$
where $b(z)$ is clearly extendible to a holomorphic function.
Here we have $\sup_{\Omega_0}|\cos z|\leq 0.913$
and, therefore, $\sup_{\Omega_0}|\cos z|+|R(z)|\leq 0.914.$
Then for $z\in\Omega_0$ we get
$$
|b(z)|=|\arccos(\cos z+R(z))-z|\leq
\int_0^1 \left|\frac{R(z)}{\sqrt{1-(\cos z+tR(z))^2}}
\right|dt
\leq 6.1 \tgo\,.
$$
\\
We now prove that 
$b(z)$ is  extendible to a holomorphic function
for $|z|<1/2$.
First we prove that there exists a real analytic positive function 
$d$ with holomorphic extension on $|z|<1/2$
such that ${w}(z)=1-\frac12\big((z-x_M)d(z)\big)^2$. By Taylor's expansion
at $z=x_M$ we have that
$
d^2(z)=-2\int_0^1(1-t){w}''(x_M+t(z-x_M))dt
$
and, therefore, for $|z|<1/2$
$$
|d^2(z)-1|\leq 1-\cos x_M +\sup_{|\z|<1/2}|\sin \z||z-x_M|+2\tgo
\leq 0.55\,.
$$
 Then we can take the principle square root\footnote{Namely taking a cut
in the negative real line.} of $d^2(z)$ obtaining the function $d(z)$.
Now consider the holomorphic function $a(z)$ define for $|z|<2$
 such that 
 $a'(z)=1/\sqrt{1-(z/2)^2}$ and $a(0)=0$.
 Then for $x$ real
  we get $a(x)=\sign(x)\arccos(1-x^2/2)$
 and also (being $d(x)>0$)
$$
b(x):=\sign(x-x_M)\arccos(1-\frac12\big((x-x_M)d(x)\big)^2)-x
=a\big((x-x_M)d(x)\big)-x\,.
$$
Then $a\big((z-x_M)d(z)\big)-z$ is a holomorphic extension 
of $b$ for $|z|<2$. An analogous argument holds for 
$|z-\pi|<2.$ 
\\
In the following we will estimates $b(z)$ for a strip
$|z|<1/2$, analogous arguments holds for $|z-\pi|<1/2$.
We will often use that\footnote{Using that $\frac12 |z|^2- (\cosh |z|-1-\frac12 |z|^2)\leq |1-\cos z|\leq \cosh |z|-1$} 
\begin{equation}\label{acqua}
|z|\leq 1\quad\Longrightarrow\quad
0.45 |z|^2\leq |1-\cos z|\leq 0.55 |z|^2
\end{equation}
Now we prove that there exists a unique function $b(z)$ defined 
for 
$$\Omega_1:=\{3\sqrt\tgo<|z|<1/2\}$$
satisfying $\sup_ {\Omega_1}|b|\leq \frac32\sqrt\tgo,$
 such that ${w}(z)=\cos(z+b(z))$, as a fixed point of the equation
 $$
 b(z)=\Psi(b)(z):=2\arcsin\left(\frac{-R(z)}{2\sin(z+b(z)/2)}\right)\,.
 $$
 Indeed,
 $$\ts
 \cos(z+b(z))-\cos z=-2 \sin\left(z+b(z)/2\right)\sin\left(b(z)/2\right)
 =R(z)\,.
 $$
 For $z\in$ we have $|z+b(z)/2|\geq \frac32\sqrt\tgo$, which 
 implies\footnote{Using that for $|z|<1$ we have $\frac45|z|\leq |\sin z|\leq\frac65 |z|$.}
 $|\sin(z+b(z)/2)|\geq \frac65\sqrt\tgo$
 and\footnote{Using that for $|z|\leq1/2$ we have $|\arcsin z|\leq \frac{2}{\sqrt 3}|z|.$}
 $
 \sup_ {\Omega_1}|\Psi(b)(z)|<\sqrt\tgo\,.
 $
 Finally $\Psi$ is a contraction 
 since\footnote{Using that for $|z|\leq1$ we have $|\cos z|\leq \sqrt 3.$}
 \begin{eqnarray*}
 &&\ts
 \sup_{\Omega_1}|\Psi(b)-\Psi(b')|
 \leq \frac{2}{\sqrt 3}\tgo
 \sup_{\Omega_1}
 \left|\frac{1}{\sin(z+b(z)/2)}-\frac{1}{\sin(z+b'(z)/2)}\right|
 \\
 &&\ts
 \leq 
 \frac{2}{\sqrt 3}\frac{5^2}{6^2}2 \left|
 \sin\left(\frac{b'(z)-b(z)}{4}\right)
 \cos\left(z+\frac{b'(z)+b(z)}{4}\right)
 \right|
 \leq\frac56 \sup_ {\Omega_1}|b-b'|\,.
\end{eqnarray*}
In conclusion we get $\sup_ {\Omega_1}|b|\leq \frac32\sqrt\tgo.$

\nl
Next, we claim that in the domain $\Omega_2:=\{|z|\leq 3\sqrt\tgo\}$ we have that 
$|b(z)|< 9\sqrt\tgo$. Indeed, by contradiction, assume that there exists $z_0\in\Omega_2$ such that for every $|z|<|z_0|$  we have $|b(z)|< 9\sqrt\tgo$ but $|b(z_0)|= 9\sqrt\tgo$.
Then $|z_0+b(z_0)|\leq 12\sqrt\tgo$ and by \eqref{acqua}
and since $\cos(z_0+b(z_0))-1=\cos z_0-1+R(z_0)$
we get
\begin{eqnarray*}
&&16\tgo\leq
0.45(|b(z_0)|-|z_0|)^2
\leq
0.45|z_0+b(z_0)|^2
\leq
|\cos(z_0+b(z_0))-1|
\\
&&\leq |\cos z_0-1|+|R(z_0)|\leq 0.55|z_0|^2+\tgo\leq 6 \tgo\,,
\end{eqnarray*}
which is a contradiction. Thus $\sup_{\Omega_2}|b(z)|\leq 9\sqrt\tgo.$ \qed

\footnotesize

\nl
{\bf Acknowledgements} The authors are grateful to  A. Neishtadt for providing parts of his Thesis (\cite{Nei89},  in Russian) related to the present paper.

\end{document}